\documentclass[11pt,reqno]{amsart}
\usepackage{amsmath,amsfonts,amssymb,amscd,amsthm,amsbsy}
\usepackage[dvips]{graphicx}
\usepackage{mathabx} \usepackage{comment} \usepackage[dvips]{color}
\usepackage[dvipsnames]{xcolor} \usepackage{fancyhdr}
\usepackage{hyperref}
\usepackage{calrsfs}
\DeclareMathAlphabet{\pazocal}{OMS}{zplm}{m}{n}

\newtheorem{theorem}{Theorem}[section]
\newtheorem{lemma}[theorem]{Lemma}

\newtheorem{proposition}[theorem]{Proposition}
\newtheorem{remark}[theorem]{Remark}
\newtheorem{definition}[theorem]{Definition}

\newcommand{\red}{\textcolor{red}}

\newcommand{\R}{\mathbb{R}}
\newcommand{\N}{\mathbb{N}}

\newcommand{\Z}{\mathbb{Z}}
\newcommand{\T}{\mathbb{T}}
\def\ep{\varepsilon}
\def\k0{K_0}
\def\kh{\widehat K}
\def\om{\omega}
\def\omh{\widehat \omega}
\def\th{\theta}
\def\kt{\widetilde K}
\def\ga{\gamma}
\def\pt{\Pi^{\top}}
\def\pp{\Pi^{\bot}}
\def\gaep{\Gamma^\ep}
\def\P{\mathcal{P}}
\def\B{\mathcal{B}}
\DeclareMathOperator{\Lip}{Lip}
\def\Span{\text{Span}}
\def\D{\pazocal{D}}

\begin{document}
\title[Persistence of Periodic Orbits] {Persistence and Smooth
  Dependence on Parameters of Periodic Orbits in Functional
  Differential Equations Close to an ODE or an Evolutionary PDE}

\author[J. Yang]{Jiaqi Yang}
\address{
School of Mathematics,
Georgia Institute of Technology,
686 Cherry St. Atlanta GA. 30332-0160 }
\email{jyang373@gatech.edu}

\author[J. Gimeno]{Joan Gimeno} 
\address{Department of Mathematics, University of Rome Tor Vergata,
  Via della Ricerca Scientifica 1, 00133 Rome (Italy)}
\email{joan@maia.ub.es}

\author[R. de la Llave]{Rafael de la Llave}
\address{
School of Mathematics,
Georgia Institute of Technology,
686 Cherry St. Atlanta GA. 30332-0160 }
\email{rafael.delallave@math.gatech.edu}

\begin{abstract}
We consider functional differential equations(FDEs) which are
perturbations of smooth ordinary differential equations(ODEs). The FDE
can involve multiple state-dependent delays or distributed delays
(forward or backward). We show that, under some mild assumptions, if
the ODE has a nondegenerate periodic orbit, then the FDE has a smooth
periodic orbit. Moreover, we get smooth dependence of the periodic
orbit and its frequency on parameters with high regularity.

The result also applies to FDEs which are perturbations of some
evolutionary partial differential equations(PDEs).

The proof consists in solving functional equations satisfied by the
parameterization of the periodic orbit and the frequency using a fixed
point approach. We do not need to consider the smoothness of the
evolution or even the phase space of the FDEs.
\end{abstract}


\maketitle

\fancyhf{}

\section{Introduction}
\label{sec:intro} 
In this paper, we first present a systematic approach to the study of
periodic orbits of functional differential equations(FDEs) which are
perturbations of smooth ordinary differential equations in $\R^n$.
This is a singular perturbation problem since the phase space of the
FDEs is infinite dimensional even if the perturbation looks small.

The approach we use bypasses completely the study of the evolution of
FDEs and we do not even need to identify the phase space. In contrast
with the standard procedure of constructing all the solutions and
selecting the periodic ones, we start with the space of periodic
functions and impose that they are solutions.

We formulate functional equations satisfied by parameterizations of
the periodic orbits and their frequencies in appropriate spaces of
smooth functions. We solve the functional equation using a fixed point
approach, which gives existence of smooth solutions and dependence on
parameters with high regularity.

One advantage of our approach is that the functional perturbations we
cover can be rather general. For example, it may include multiple
state-dependent delays, distributed delays, or implicitly defined
delays of the type appearing in electrodynamics (see Section
\ref{sec:electrodynamics}). The delays can be either regular backward
delays or forward delays (as in advanced equations).

Then, using a similar but more elaborate proof, we get results on
periodic orbits for equations with small delays, which have
applications in electrodynamics.

Finally, we extend the results to perturbations of partial
differential equations(PDEs). We can consider PDEs which have good
forward (but not backward) evolutions such as parabolic equations as
well as some ill-posed equations (e.g. Boussinesq equation in water
waves, which even if ill posed, admits many physically interesting
solutions).

A philosophy similar to that of this paper has been also used in other
papers. \cite{RaXueD, HR, HR2} develop functional equations for
quasi-periodic solutions in several contexts and study them using KAM
theory.  In \cite{theo}, one can find a theory of persistence of
stable manifolds in some limited contexts.  We hope that some of the
previous studies can be extended to more dynamical objects. Notably we
expect to get higher regularity of the center manifolds for SDDEs,
which is essential for applications of the center manifold reduction
to bifurcation theory \cite{Hum}.  Of course, removing the
perturbative setting remains a long term goal, but this seems to pass
through refining the theory of existence and regularity of
\cite{Walt}.

\subsection{Backgrounds on Functional Differential Equations}
In many applications, one needs to consider FDEs. Delay differential
equations(DDEs) appear naturally as models in electrodynamics, control
theory, biology, neuroscience, and economics, see \cite {WheelerF49,
  DriverElec, Driver, hart, Hum, Campbell, CampbellZhu, Freedman,
  WalM, mack} and references therein.  In many cases, the delays
depend on the states of the systems, one needs to consider
state-dependent delay equations(SDDEs). For example, in the
formulation of Electrodynamics, the delays depend implicitly on the
solution. Sometimes several delays are involved in one equation, with
different forms. Besides the interest in applications, the field of
FDEs is a very rich mathematical subject worth of study because of its
own depth.

The theory for delay equations with constant delays is well
established \cite{HaleLunel, Diek}. However, many fundamental problems
are not settled for SDDEs. For example, even identifying the correct
phase space to formulate the equation is not clear.  The paper
\cite{Walt} made a breakthrough considering a submanifold of $C^1$
space, the solution manifold, as phase space for SDDEs on which the
semiflow is $C^1$. A result on differentiability of solutions with
respect to parameters for a class of SDDEs in Sobolev sense (using
quasi-Banach spaces) is in \cite{HartungTuri}. It seems that there is
no result on higher regularity of the semiflow and dependence on
parameters for a general solution. One can refer to \cite{hart} for a
review of the applications and results in SDDEs.  SDDEs display rich
behaviors, see \cite{Hum, HumphriesDMU12}. At the same time, some
SDDEs, like the ones considered in this paper, have many solutions
with regular behaviors, see \cite{RaXueD, HR, HR2, theo}. See also
\cite{Pablo, MCR} for results on low regularity via topological
approach.

Periodic orbits are important landmarks in dynamical systems. There
has been interest on studying periodic orbits in DDEs, see
\cite{NusPer, KapYor, MNglobal, JaqLM}. Some studies in the setting of
SDDEs are in \cite{MNAna, MPNP, MPN, SieberP}. Some numerical works
are in \cite{BGL, Sz, SzZ}.

\subsection{Related Results in the Literature}
Results on persistence of non-degenerate periodic orbits and
dependence on parameters for FDEs with constant delays was proven by
studying the evolution operator, see \cite{HaleLunel, HaleS}, and
\cite{HaleW}.  This method is difficult to apply to SDDEs for
regularities higher than $C^1$ since one would need to extend the
regularity theory of the evolution \cite{Walt} to higher regularities.
The paper \cite{Mawhin} also studied functional equations satisfied by
periodic orbits, but treats them using topological methods, which do
not allow to study regularity. See also the excellent surveys
\cite{Mawhin71,Nussbaum79}.

\subsection{Organization of the Paper}
A precise formulation of the problem is given in Section
\ref{form}. Section \ref{sec:para} introduces the parameterization
method for our problem.  Section \ref{sec:main} states the main
results of this paper.  This result is formulated in terms of
properties of the functional $P$.  The detailed proofs of the main
results are in Section \ref{sec:proof}.  In Sections
\ref{sec:examauto} and \ref{sec:examper} we verify that several models
that appear in the literature indeed satisfy the assumptions of the
main result.  These sections are the core of this paper.

In the other sections, we present extesions of the method and the
philosopy and show that they lead to results for several models in the
literature.

Section \ref{sec:small_delays} is devoted to the analysis of equations
with small delays, which requires an extension of the general result
and indeed requires stronger regularity assumptions.  Section
\ref{sec:electrodynamics} considers equations appearing in
electrodynamics, which has been a very important motivation for the
whole theory of FDEs.  In particular, we give some justification to
several procedures used in Physics such as the $1/c$ expansions.

Section \ref{sec:hyperbolic} introduces a different method for the
case that the periodic orbits are hyperbolic.  Even if this is a
particular case of the previous results for ODEs, it generalizes to
evolutionary Partial Diffferential Equations.  In Section
\ref{sec:pde}, we present results for several evolutionary PDEs which
have received attention in the literature.  We note that, since our
method dispenses with defining the evolution, the results apply even
to ill-posed PDEs.

In Appendix~\ref{app:regularity} we have collected some results of
analysis that we need to use.

\section{Formulation of the Problem}\label{form}

Consider an $n$-dimensional ODE
\begin{equation}\label{ode}
\dot x(t)=f(x(t)),
\end{equation}
where, for the moment, $f\colon\R^n\to \R^n$ is a $C^{\infty}$ vector
field (later we will assume less regularity).

We assume that equation \eqref{ode} has a periodic orbit with
frequency $\omega_0\neq 0$. The existence of periodic solutions for
ODEs will not be discussed here. (We note however that the same
methods discussed here can be used to produce periodic solutions of
the ODE perturbatively.)

We consider singular perturbation of equation \eqref{ode} to FDEs with
parameter $\ga$:
\begin{equation}\label{fde}
\dot x(t)=f(x(t))+\ep P(x_t,\ga),
\end{equation}
where $P\colon \mathcal{R}[-h,h]\times O\to \R^n$, $h$ is a positive
constant. $\mathcal{R}[-h,h]$ is a space of regular functions from
$[-h,h]$ to $\R^n$.  The precise regularity of the functions in
$\mathcal{R}[-h,h]$ will be specified later.  The \emph{``history
  segment''} $x_t\in \mathcal{R}[-h,h]$ is defined as $x_t(s)=x(t+s)$
for $s\in [-h,h]$. The parameter $\gamma\in O$, where $O$ is a bounded
open set in $\R^m$.  Note that we allow that our history segments
involve also the future, so that the theory we will develop applies
not just to delay equations but to equations that involve the future.

In many treatments of delay equations it is customary to think of
$\mathcal{R}[-h,h]$ as the phase space in which one sets initial
conditions and defines an evolution.  For example, in the case of
constant delay equations, it is customary to impose initial conditions
in $C^0[-h,0]$, with constant $h$ being the delay. Nevertheless, in
the case of SDDEs, this space includes many functions which cannot
satisfy the equations and, therefore, have no physical meaning.  As it
will be clear later, our treatment bypasses the consideration of the
evolution defined by the ODE, so that we will not think of
$\mathcal{R}[-h,h]$ as the phase space of the evolution.

Under nondegeneracy condition on the periodic orbit of equation
\eqref{ode} and some mild assumptions on $P$, see more details in the
definition of $\P$ in \eqref{eqk} and assumptions \eqref{H2.1},
\eqref{H3.1}, \eqref{H2.2}, and \eqref{H3.2}, we show that for small
enough $\ep$, there exists periodic orbit for FDE \eqref{fde}. We also
show that the periodic orbits for equation \eqref{fde} depend on $\ga$
smoothly.

From now on, we will identify the periodic orbit in a function space
for FDE \eqref{fde} with a periodic function having values in
$\R^n$. Under this identification, we will see that the periodic orbit
for FDE \eqref{fde} is close to the periodic orbit for equation
\eqref{ode} for small $\ep$.

\section{Parameterization Method}\label{sec:para}
Let $\k0\colon \T\to \R^n$ be a parameterization of the periodic orbit
of equation \eqref{ode}, where $\T=\R/\Z$. This means that for any
fixed $\theta$, $x(t)=\k0(\theta+\omega_0t)$ solves equation
\eqref{ode}. Equivalently, $\k0$ satisfies the functional equation
(invariance equation):
\begin{equation}\label{eqk0}
\omega_0D\k0(\theta)=f(\k0(\theta)).
\end{equation}
Note that such $\k0$ is unique up to a phase shift. In this case,
$\k0$ is $C^\infty$ since $f$ is $C^\infty$.

We aim to find $K\colon\T\to \R^n$ and $\om>0$, such that for any
$\th$, $x(t)=K(\th+\om t)$ solves equation \eqref{fde}. And we say
such $K$ parameterizes the periodic orbit of FDE \eqref{fde}.

$x(t)=K(\th+\om t)$ solving equation \eqref{fde} is equivalent to $K$
satisfying the functional equation:
\begin{equation}\label{eqk}
\om D K(\th)=f(K(\th))+\ep \P(K, \om, \ga, \th),
\end{equation}
$\P(K, \om, \ga,\th)$ results from substituting $x(t)=K(\th+\om t)$
into $P(x_t,\ga)$ in equation \eqref{fde} and letting $t=0$. See
Sections \ref{sec:examauto}, \ref{sec:examper}, and
\ref{sec:small_delays} for explicit formulations of $\P$ in some
specific examples.

The equation \eqref{eqk} will be the centerpiece of our treatment.  We
will see that, using different methods of analysis, we can give
results on existence of solutions of \eqref{eqk}.  Note that this
analysis produces periodic solutions of \eqref{fde} without discussing
a general theory of existence and dependence on parameters of the
solutions.
\section{Main Results}\label{sec:main}

\subsection{Assumptions}\label{sec:assump}
For a given $\th_0\in \T$, let $\Phi(\theta;\theta_0)$ be the
fundamental solution of the variational equation of the ODE
\eqref{ode}, i.e.,
\begin{equation}\label{var}
\omega_0\frac{d}{d\theta}\Phi(\theta;\theta_0) =
Df(\k0(\theta))\Phi(\theta;\theta_0),\quad \Phi(\theta_0;\theta_0)=Id.
\end{equation}

We need to assume that the periodic orbit of \eqref{ode} is
nondegenerate, that is we impose the following assumption on
$\Phi(\theta_0+1;\theta_0)$:
\smallskip
\begin{enumerate}
\renewcommand{\theenumi}{H\arabic{enumi}}
\renewcommand{\labelenumi}{(\theenumi)}
\item \label{H1} $\Phi(\theta_0+1;\theta_0)$ has a simple eigenvalue 1
  whose eigenspace is generated by $DK_0(\th_0)$.
\end{enumerate}

Note that, because of the existence and uniqueness of the solutions of
\eqref{var}, and the periodicity of $\k0$, we have that
\[
\begin{split} 
 &\Phi(\theta_2; \theta_0) = \Phi(\theta_2; \theta_1) \Phi(\theta_1;
  \theta_0); \\ & \Phi(\theta_1 + 1 ; \theta_0 + 1 ) = \Phi(\theta_1;
  \theta_0). \\
\end{split} 
\] 

As a consequence, 
\[
\Phi(\theta_0 + 1 ; \theta_0) = \Phi(\theta_0 +1; 1) \Phi(1 ; 0 )
\Phi(0 ; \theta_0) = \Phi(0;\theta_0)^{-1} \Phi(1 ; 0 ) \Phi(0 ;
\theta_0) .
\] 
So that the spectrum of $\Phi(\theta_0 +1; \theta_0)$, commonly called
the Floquet multipliers, is independent of the starting point
$\theta_0$.

Under assumption \eqref{H1}, there exists an $(n-1)$-dimensional
linear space $E_{\th_0}$ at $\k0(\th_0)$, \big(the spectral complement
of $\Span\{ D\k0(\th_0)\}$, corresponding to the eigenvalues of
$\Phi(\theta_0+1;\theta_0)$ other than $1$, $\R^n=E_{\th_0}\oplus
\Span\{ D\k0(\th_0)\}$\big), on which $\left[ Id-\Phi(\th_0+1;\th_0)
  \right]$ is invertible. We denote the the projections onto $\Span\{
D\k0(\th_0)\}$ and $E$ as $\pt_{\th_0}$ and $\pp_{\th_0}$,
respectively.

\begin{remark}
An equivalent formulation of \eqref{H1} in terms of functional
analysis is \eqref{H1'}. Define the operator $\mathcal{L}\colon
C^1(\T,\R^n)\to C^0(\T,\R^n)$:
\[
\mathcal{L}(v)(\theta)=\omega_0Dv(\theta)-Df(\k0(\theta))v(\theta).
\]
\begin{enumerate}
\renewcommand{\theenumi}{H\arabic{enumi}'}
\renewcommand{\labelenumi}{(\theenumi)}
\item \label{H1'} $Range(\mathcal{L})$ is co-dimension 1,
  $Range(\mathcal{L})\oplus \Span\{D\k0\}=C^0(\T,\R^n)$.
\end{enumerate}
The proofs of the Theorems in the next section imply the equivalence
of \eqref{H1} and \eqref{H1'}.
\end{remark}

\medskip
 
To show the persistence of periodic orbit for a fixed $\ga\in O$, the
following assumptions on $\P$ are crucial.  The assumption
\eqref{H2.1} is about smoothness of $\P$ and expresses that $\P$ maps
$C^{\ell+\Lip}$ balls around zero into $C^{\ell-1+\Lip}$ balls around
zero. \eqref{H3.1} is about Lipschitz property of $\P$ in $C^0$ for
smooth $K$'s. These properties are verified in the examples we study
in Sections \ref{sec:examauto} and \ref{sec:examper}. For example,
when the functional $P$ is evaluation on $x(t-r(x(t)))$, the
regularity is a consequence of the fact that we can control the
$C^\ell $ norm of $f\circ g$ by the $C^\ell$ norm of $f,g$. (We can
even loose a derivative).  The $C^0$ Lipschitz property results from
the mean value theorem ($\|f\circ g_1 - f\circ g_2 \|_{C^0} \le
\|f\|_{C^1}\|g_1 - g_2\|_{C^0}$).

In the following, $\ell$ is an arbitrarily fixed positive integer. 

Let $U_\rho$ be the ball of radius $\rho$ in the space
$C^{\ell+\Lip}(\T, \R^n)$ centered at $\k0$, and let $B_\delta$ be the
interval in $\R$ with radius $\delta$ centered at $\om_0$.

\begin{enumerate}
\setcounter{enumi}{1}
\renewcommand{\theenumi}{H\arabic{enumi}.1}
\renewcommand{\labelenumi}{(\theenumi)}
\item \label{H2.1} If $K\in U_\rho$ and $\om\in B_\delta$, then $\P(K,
  \om, \ga, \cdot)\colon \T\to \R^n$ is $C^{\ell-1+\Lip}$, with
  $\|\P(K, \om, \ga, \cdot)\|_{C^{\ell-1+\Lip}}\leq
  \phi_{\rho,\delta}$, where $\phi_{\rho,\delta}$ is a positive
  constant that may depend on $\rho$ and $\delta$.

\item \label{H3.1} For $K,~K'\in U_\rho$, and $\om,~\om' \in
  B_\delta$, there exists constant $\alpha_{\rho,\delta}>0$, such that
  for all $\th\in\T$,
\begin{equation}
\phantom{AAA}|\P(K, \om, \ga, \th)-\P(K', \om', \ga, \th)|\leq
\alpha_{\rho,\delta}\max\left\{|\om-\om'|,\|K-K'\|\right\},
\end{equation}
where $\|K-K'\|$ is the $C^0$-norm of $K-K'$ under the Euclidean
distance on $\R^n$.
\end{enumerate}

To show that the periodic orbits of the FDE \eqref{fde} depend on
parameter $\ga$ smoothly, one needs to consider $K$ as a function of
$\th$ and $\ga$, and $\om$ as a function of $\ga$. \eqref{H2.2} and
\eqref{H3.2} are similar to \eqref{H2.1} and \eqref{H3.1},
respectively.

We let $\pazocal{U}_\rho$ be the ball of radius $\rho$ in the space
$C^{\ell+\Lip}(\T\times O, \R^n)$ centered at $\k0$, and let
$\pazocal{B}_\delta$ be the ball in $C^{\ell+\Lip}(O, \R)$ with radius
$\delta$ centered at constant function $\om_0$.
\begin{enumerate}
\setcounter{enumi}{1}
\renewcommand{\theenumi}{H\arabic{enumi}.2}
\renewcommand{\labelenumi}{(\theenumi)}
\item \label{H2.2} If $K\in\pazocal{U}_\rho$ and
  $\om\in\pazocal{B}_\delta$, then $\P(K, \om, \cdot, \cdot)\colon
  \T\times O\to \R^n$ is $C^{\ell+\Lip}$ in $\ga$, and
  $C^{\ell-1+\Lip}$ in $\th$, with the bounds $\|\P(K, \om, \cdot,
  \th)\|_{C^{\ell+\Lip}}\leq \phi_{\rho,\delta}$, and $\|\P(K, \om,
  \ga, \cdot)\|_{C^{\ell-1+\Lip}}~\leq \phi_{\rho,\delta}$, where
  $\phi_{\rho,\delta}$ is a positive constant.
\item \label{H3.2} For $K, ~K'\in\pazocal{U}_\rho$ and $\om,
  ~\om'\in\pazocal{B}_\delta$, there exists constant
  $\alpha_{\rho,\delta}>0$, such that for all $\th\in\T$ and $\ga\in
  O$,
\[
\phantom{AAAA}|\P(K, \om, \ga, \th)-\P(K', \om', \ga, \th)|\leq
\alpha_{\rho,\delta}\max\left\{\|\om-\om'\|,\|K-K'\|\right\},
\]
where $\|\om-\om'\|$ is the $C^0$-norm of $\om-\om'$.
\end{enumerate}

\begin{remark}\label{ep}
Note that our results work exactly the same if the perturbation
depends on $\ep$, i.e. we have $P(x_t, \ga, \ep)$ instead of
$P(x_t,\ga)$ in \eqref{fde}. We can get $\P(K, \om, \ga, \ep, \th)$ in
this case. We need assumptions on $\P$ to hold uniformly in $\ep$ for
all small $\ep$.

More specifically, \eqref{H2.1}, \eqref{H3.1} can be reformulated as:
\begin{enumerate}
\setcounter{enumi}{1}
\renewcommand{\theenumi}{H\arabic{enumi}.1'}
\renewcommand{\labelenumi}{(\theenumi)}
\item \label{H2.1'} If $K\in U_\rho$ and $\om\in B_\delta$, then
  $\P(K, \om, \ga, \ep, \cdot)\colon \T\to \R^n$ is $C^{\ell-1+\Lip}$,
  with \[\|\P(K, \om, \ga, \ep, \cdot)\|_{C^{\ell-1+\Lip}}\leq
  \phi_{\rho,\delta}(\ep).\] Function $\phi_{\rho,\delta}$ satisfies
  that $\ep\phi_{\rho,\delta}(\ep)$ converges to zero as $\ep\to 0$.
\item \label{H3.1'} For $K,~K'\in U_\rho$, and $\om,~\om' \in
  B_\delta$, there exists positive function $ \alpha_{\rho,\delta}$,
  such that for all $\th\in\T$,
\[
\phantom{AAAA} |\P(K, \om, \ga, \ep, \th)-\P(K', \om', \ga, \ep,
\th)|\leq
\alpha_{\rho,\delta}(\ep)\max\left\{|\om-\om'|,\|K-K'\|\right\},
\]
function $ \alpha_{\rho,\delta}$ satisfies that $\ep
\alpha_{\rho,\delta}(\ep)$ converges to zero as $\ep\to 0$.
\end{enumerate}
The assumptions similar to \eqref{H2.2}, \eqref{H3.2} can be
formulated similarly.

\end{remark}

\begin{remark} 
The assumptions we use are similar to assumptions in invariant
manifold theory.  For example in \cite{Lan}, the \eqref{H2.1} is
called \emph{propagated bounds}.
\end{remark}

\begin{remark}
We call attention to the fact that in Section~\ref{sec:small_delays}
we will weaken substantially the assumption \eqref{H3.1} to be able to
deal with equations with small delays.

\end{remark} 

\subsection{Main Theorems}

Let $\N$ denote the set of positive numbers.
\begin{theorem}[Persistence]\label{thm:per}
For a given $\ell\in \N$, assume that $f$ in \eqref{fde} is
$C^{\ell+\Lip}$, and that \eqref{H1}, \eqref{H2.1}, and \eqref{H3.1}
are satisfied for a given $\ga\in O$. Then, there exists $\ep_0>0$,
such that when $\ep<\ep_0$, the FDE \eqref{fde} has a periodic orbit,
which is parameterized by $K\colon\T\to \R^n$. The smallness condition
of $\ep_0$ depends on $\ell$, $f$, and $\P$.

The frequency $\om$ for the periodic orbit is close to $\om_0$, the
frequency of the periodic orbit of equation \eqref{ode}.
$\|K-\k0\|_{C^{\ell}}$ is small under a suitable choice of the phases.

\end{theorem}

\begin{theorem}[Smooth Dependence on Parameter]\label{thm:smooth}
For a given $\ell\in \N$, assume that $f$ in \eqref{fde} is
$C^{\ell+\Lip}$, and that \eqref{H1}, \eqref{H2.2}, and \eqref{H3.2}
are satisfied. Then, there is $\ep_0>0$, such that if $\ep<\ep_0$, one
can find $K_{\gamma}(\theta)$ which parameterizes the periodic orbit
of FDE \eqref{fde} persisted from the periodic orbit of
\eqref{ode}. The smallness condition of $\ep$ depends on $\ell$, $f$,
and $\P$.

$K_{\gamma}$ has frequency $\omega_{\gamma}$. $K_{\gamma}(\theta)$ is
jointly $C^{\ell+\Lip}$ in $\theta$ and $\gamma$, $\omega_\ga$ is
$C^{\ell+\Lip}$ in $\gamma$.
\end{theorem}

\subsection{Some Comments on the Theorems \ref{thm:per} and \ref{thm:smooth}}

\begin{remark}\label{rmk:energy} 
One physically important case where assumption \eqref{H1} fails is
when there is a conserved quantity (for example, the energy in
mechanical systems). We are not able to deal with this case by the
method of this paper, but we hope to come back to this problem.
\end{remark}

\begin{remark}\label{rmk:phase}
Note that $K$ will not be unique. If $K(\th)$ parameterizes the
periodic orbit, then for any given $\th_1$, $K(\th+\th_1)$ also
parameterizes the periodic orbit, with a shifted phase. Hence, in
Theorem \ref{thm:per}, the smallness of $\|K-\k0\|_{C^{\ell-1}}$ is
interpreted under a suitable choice of the phases.

This is the only source of non-uniqueness since the proofs of Theorems
\ref{thm:per} and \ref{thm:smooth} are based on contraction mapping
argument, the parameterizations we found are locally unique up to
phase shifts.
\end{remark}

\begin{remark}
The smallness of $\ep$ depends on $\ell$, hence, the method cannot get
a $C^\infty$ result directly. Note, however, that in some cases,
e.g. state-dependent delay perturbations in equation \eqref{sdde}, one
can bootstrap the regularity from $C^1$ to $C^\infty$.
\end{remark}

\begin{remark}\label{rmk:gen}
Our results apply to several types of FDEs, especially to many DDEs,
see Sections \ref{sec:examauto} and \ref{sec:examper}. We only need
that \eqref{H1}, \eqref{H2.1}, \eqref{H2.2}, \eqref{H3.1}, and
\eqref{H3.2} are satisfied. Indeed, we allow several terms in the
equation which may involve forward and backward delays.
\end{remark}

\begin{remark}
Our method allows to bypass the propagation of discontinuity in
DDEs. Moreover, it has no restriction on the relation between the
period of the periodic orbits and the size of the delay.
\end{remark}

\begin{remark}
The proofs we present are constructive, hence they can be implemented
numerically. Indeed, we formulate the problem as a fixed point of a
contractive operator, which concatenates several elementary
operators. Implementations of these elementary operators for a 2D
model are addressed in a numerical toolkit developed in \cite{Joan19}.

The proofs, based on fixed point approach, also lead to results in an
a-posteriori format, which state that if there is an approximate
solution (satisfying some mild assumptions), then there is a true
solution which is close to the approximate one. See more details in
Section \ref{sec:conc}.
\end{remark}

\begin{remark}
A-posteriori theorems justify asymptotic expansions where solutions
are written as formal expansions in terms of the small parameters, see
\cite{Chiconedelay, CasalCL20}. Truncations of the formal power series
provide approximate solutions. The a-posteriori theorem shows that
there is one true solution close by.

A-posteriori theorems are also the base of computer-assisted proofs.
Numerical methods produce approximate solutions. If one can estimate
rigorously the error and the non-degeneracy conditions, then one has
established the existence of the solution. The verification of the
error in the approximation is a finite (but long) calculation which
can be done using computers taking care of round-off and
truncation. Some cases where computer assisted proofs have been used
in constant delay equations for periodic orbits and unstable manifolds
are \cite{KissL12, GroothedeMJ17}.
\end{remark} 

\section{Proofs}\label{sec:proof}

The proof of Theorems \ref{thm:per} and \ref{thm:smooth} are based on
fixed point approach. We will provide the detailed proof of Theorem
\ref{thm:per}. The proof of Theorem \ref{thm:smooth} follows in the
same manner by adding the parameters in the unknowns, see
\ref{sec:cosm}.

The proof consist of several steps. First, we define an operator in an
appropriate space of smooth functions. Then, we show that (i) the
operator maps a ball in this space into itself (Section
\ref{sec:prop}); (ii) the operator is a contraction in a $C^0$ type of
distance (Section \ref{sec:cont}). The existence of fixed point in
desired space is hence ensured using a generalization of contraction
mapping \cite{Lan}.

\subsection{Invariance Equations}\label{sec:inveq}
In this section, we reformulate the invariance equation \eqref{eqk}.
Since we expect that the solutions $K, \omega$ will be small
perturbations of the unperturbed ones, it is natural to reformulate
\eqref{eqk} as an equation for the corrections from the unperturbed
ones. In Section ~\ref{sec:operator} we will manipulate the equation
for the corrections into a fixed point problem.

Let
\begin{equation}\label{kom}
\begin{split}
K(\theta)&\coloneq \k0(\theta)+\kh(\theta),\\ \om&\coloneq\om_0+\omh,
\end{split}
\end{equation}
where $\kh\colon\T\to \R^n$ and $\omh\in\R$ are corrections to the
parameterization and frequency of the periodic orbit of the
unperturbed equation. Our goal is to find $\kh$ and $\omh$ so that $K$
and $\om$ satisfy the functional equation \eqref{eqk}.

Using the notation in \eqref{kom} and the invariance equation
\eqref{eqk0} for $\k0$ and $ \om_0$, we are led to the following
functional equation for $\kh$ and $\omh$,
\begin{equation}\label{inveq}
\om_0 D\kh(\th)-Df(\k0(\th))\kh(\th)=B^\ep(\kh,\omh,\ga,\th)-\omh D\k0(\th),
\end{equation}
where
\begin{align}\label{Bdef}
&B^\ep(\kh,\omh,\ga,\th)\coloneq N(\th, \kh)+\ep\P(K,\om,\ga,\th)-\omh
  D\kh(\th),\\ &N(\th, \kh)\coloneq
  f(\k0(\th)+\kh(\th))-f(\k0(\th))-Df(\k0(\th))\kh(\th).\nonumber
\end{align}

The basic idea for this regrouping is that since $K$ and $\om$ are
expected to be close to $\k0$ and $\om_0$, we only need to find the
corrections.

\subsection{The Operator}
\label{sec:operator}
Recall $\Phi(\theta;\theta_0)$ introduced in \eqref{var} as the flow
of the variational equations. Using the variation of parameters
formula, equation \eqref{inveq} for $\kh$ and $\omh$ is equivalent to:
\begin{equation}\label{fixedpteq}
\kh(\th)=\Phi(\th;\th_0)\left\{u_0+\frac{1}{\om_0}\int^{\th}_{\th_0}\Phi(s;\th_0)^{-1}\big(B^\ep(\kh,\omh,\ga,s)-\omh
D\k0(s)\big)ds\right\},
\end{equation}
where the initial condition $\kh(\th_0)=u_0$ is to be found imposing
that $\kh$ is periodic.  This will be discussed in the Section
\ref{sec:periodiccon}.

We can think of \eqref{fixedpteq} as a fixed point equation. The right
hand side is an operator in $\kh$, see Section \ref{definitionop}. We
start with a given $\kh$, choose $\omh$ following Section
\ref{sec:periodiccon} and we substitute them in right hand side of
\eqref{fixedpteq}.

\subsubsection{Periodicity Condition}\label{sec:periodiccon}

Since the right hand side of equation \eqref{inveq} is periodic, $\kh$
is periodic if and only if $\kh(\th_0)=\kh(\th_0+1)$, i.e.,
\begin{align}\label{per1}
[Id-\Phi(\th_0+1;\th_0)]u_0=\frac{1}{\om_0}&\Phi(\th_0+1;\th_0)\int^{\th_0+1}_{\th_0}\Phi(s;\th_0)^{-1}B^\ep(\kh,\omh,\ga,s)ds\nonumber\\ -&\frac{\omh}{\om_0}\Phi(\th_0+1;\th_0)\int^{\th_0+1}_{\th_0}\Phi(s;\th_0)^{-1}
D\k0(s)ds.
\end{align}

Since $\k0$ solves \eqref{eqk0}, we have
\[
\Phi(s;\th_0)D\k0(\th_0)= D\k0(s).
\]

Then, the periodicity condition \eqref{per1} becomes
\begin{align}\label{per2}
[Id-\Phi(\th_0+1;\th_0)]u_0=\frac{1}{\om_0}&\int^{\th_0+1}_{\th_0}\Phi(\th_0+1;s)B^\ep(\kh,\omh,\ga,s)ds\nonumber\\ -&\frac{\omh}{\om_0}
D\k0(\th_0).
\end{align}

One is able to solve for $u_0$ if the right hand side of equation
\eqref{per2} is in the range of $Id-\Phi(\th_0+1;\th_0)$. Thanks to
assumption \eqref{H1}, this can be achieved by choosing $\omh$. Such
$\omh$ is unique.

\subsubsection{Spaces}
\label{sec:defspaces}
Let $I_a=[-a, a]$ be an interval which contains $0$, where $a>0$, let 
\begin{equation}\label{space}
 \begin{split}
  \mathcal{B}_\beta=\Big\{g\colon \T\to\R^n\big|~& \text{$g$ is
    $C^{\ell+\Lip}$}, ~\Big\|\frac{d^i}{d\th^i}g(\th)\Big\|\leq
  \beta_i,\ i=0,1,\dotsc,\ell,\ \\ &
  \Lip\Big(\frac{d^l}{d\th^l}g(\th)\Big)\leq \beta^{\Lip}_\ell\Big\},
 \end{split}
\end{equation}
where $\beta=(\beta_0, \beta_1,\dotsc, \beta_\ell,
\beta^{\Lip}_\ell)$. The constants $a$, $\beta_i$,
$i=0,1,\dotsc,\ell$, and $\beta^{\Lip}_\ell$ will be chosen in the
proof.
\subsubsection{Definition of the Operator}
\label{definitionop}
Define the operator $\Gamma^\ep$ on $I_a\times \mathcal{B}_\beta$,
\begin{equation}\label{op}
\Gamma^\ep(\omh,\kh)=\begin{pmatrix}
\Gamma^\ep_1(\omh,\kh)\\ \Gamma^\ep_2(\omh,\kh)
\end{pmatrix}.
\end{equation}
Componentwise,
\begin{equation}\label{op1}
 \Gamma^\ep_1(\omh,\kh)=\frac{\langle\int^{\th_0+1}_{\th_0}\pt_{\th_0}\Phi(\th_0+1;s)B^\ep(\kh,\omh,\ga,s)ds,
   D\k0(\th_0)\rangle}{\left| D\k0(\th_0)\right|^2},
\end{equation}
where $\langle\cdot,\cdot\rangle$ is the standard inner product on
$\R^n$.

\begin{align}\label{op2}
 \Gamma^\ep_2(\omh,\kh)(\th)=\Phi&(\th;\th_0)u_0\\ +&\frac{1}{\om_0}\int^{\th}_{\th_0}\Phi(\th;s)\big(B^\ep(\kh,\omh,\ga,s)-\Gamma^\ep_1(\omh,\kh)
 D\k0(s)\big)ds\nonumber,
\end{align}
where $u_0\in E$ satisfies 
\begin{align}\label{u0}
[Id-\Phi(\th_0+1;\th_0)]u_0=&\frac{1}{\om_0}\int^{\th_0+1}_{\th_0}\Phi(\th_0+1;s)B^\ep(\kh,\omh,\ga,s)ds\nonumber\\ &-\frac{\gaep_1(\omh,\kh)}{\om_0}
D\k0(\th_0)\\ =&\frac{1}{\om_0}\int^{\th_0+1}_{\th_0}\pp_{\th_0}\Phi(\th_0+1;s)B^\ep(\kh,\omh,\ga,s)ds.\nonumber
\end{align}

\begin{remark}
The choice of $ \Gamma^\ep_1(\omh,\kh)$ ensures that the right hand
side of \eqref{u0} is in the range of $Id-\Phi(\th_0+1;\th_0)$. Since
the kernel of $Id-\Phi(\th_0+1;\th_0)$ is $\Span\{ D\k0(\th_0)\}$,
equation \eqref{u0} has infinitely many solutions, all of them are the
same up to constant multiples of $D\k0(\th_0)$. In the definition of
the operator $\gaep$, we have chosen the solution for equation
\eqref{u0} which lies in the space $E$. If we choose a different $u_0$
solving \eqref{u0}, we will get another parameterization of the
periodic orbit corresponding to a different phase, see Remark
\ref{rmk:phase}.
\end{remark}

Our goal is to find the fixed point $(\omh^*, \kh^*)$ of the operator
$\gaep$ in a ball $I_a\times\B_\beta$, which will solve the equation
\eqref{inveq}. Hence $\om=\om_0+\omh^*$ and $K=\k0+\kh^*$ satisfy
\eqref{eqk}, $K$ parameterizes the periodic orbit of \eqref{fde} with
frequency $\om$.

To this end, under the assumptions \eqref{H1}, \eqref{H2.1} and
\eqref{H2.2}, we show in Section~\ref{sec:prop} that for small $\ep$,
we can choose $a$ and $\beta$ so that $\gaep$ maps $I_a\times\B_\beta$
back into itself.

In Section~\ref{sec:cont} we show that $\gaep$ is a contraction in a
$C^0$-like distance. The desired result of existence of a locally
unique fixed point follows from a fixed point result in the literature
that we have collected as Lemma~\ref{closure}.

\subsection{Propagated Bounds for $\gaep$}\label{sec:prop}

In this section, we will prove the following Lemma.

\begin{lemma}\label{lem:prop}
Assume $\ep$ is small enough, then $a$ and $\beta$ can be chosen such
that $\gaep\colon I_a\times\mathcal{B}_\beta\to
I_a\times\mathcal{B}_\beta $.
\end{lemma}
\begin{proof}
Note that
\[
\|N(\th,\kh)\|\leq \frac{1}{2} \Lip(Df) \|\kh\|^2,
\]
where $\|\cdot\|$ means $C^0$-norm. Indeed, here and later in this
proof we only need the Lipschitz constant of $Df(x)$ in a neighborhood
of the periodic orbit of the unperturbed ODE, i.e. $\k0(\T)$.

Using the  integration by parts formula, for $\th\in[\th_0,\th_0+1]$, we have
\[
\left|\int^\th_{\th_0}\Phi(\th_0+1;s)\omh
D\kh(s)ds\right|\leq\Big(2\|\Phi(\th_0+1;\th)\|+\Big\|\frac{d}{d\th}\Phi(\th_0+1;\th)\Big\|\Big)|\omh|\|\kh\|.
\]
where 
\[
\|\Phi(\th_0+1;\th)\|\coloneq \max_{\th\in [\th_0,\th_0+1]}|\Phi(\th_0+1;\th)|,
\]
and \[ \left\|\frac{d}{d\th}\Phi(\th_0+1;\th)\right\|\coloneq
\max_{\th\in
  [\th_0,\th_0+1]}\left|\frac{d}{d\th}\Phi(\th_0+1;\th)\right|,
\]
$|\cdot|$ denotes the matrix norm. We will use similar conventions for
norms from now on.

Since $(\omh,\kh)\in I_a \times\B_\beta$, we have
\begin{align}\label{estga1}
|\gaep_1(\omh,\kh)|\leq \frac{\|\pt_{\th_0}\|}{\big| D\k0(\th_0)\big|}\Big[\|\Phi&(\th_0+1;s)\|\Big(\frac{1}{2}\Lip(Df) \beta_0^2+\ep\|\P(K,\om,\ga,\th)\|\Big)\nonumber\\
&+\Big(2\|\Phi(\th_0+1;\th)\|+\Big\|\frac{d}{d\th}\Phi(\th_0+1;\th)\Big\|\Big)a\beta_0\Big],
\end{align}
and,
\begin{align}\label{estga2}
\|\gaep_2(\omh,\kh)\|\leq \|\Phi&(\th;\th_0)\|M\Big[\|\Phi(\th_0+1;s)\|(\frac{1}{2}\Lip(Df) \beta_0^2+\ep\|\P(K,\om,\ga,\th)\|)\nonumber\\
&\phantom{AAAAAA}+\Big(2\|\Phi(\th_0+1;\th)\|+\Big\|\frac{d}{d\th}\Phi(\th_0+1;\th)\Big\|\Big)a\beta_0\Big]\nonumber\\
&+\frac{1}{\om_0}\Big[\|\Phi(\th;s)\|(\frac{1}{2}\Lip(Df) \beta_0^2+\ep\|\P(K,\om,\ga,\th)\|)\\
&\phantom{AAAAAA}+\Big(2\|\Phi(\th;\th_0)\|+\Big\|\frac{d}{ds}\Phi(\th;s)\Big\|\Big)a\beta_0\Big]\nonumber\\
&+\frac{\|\Phi(\th;s)\|}{\om_0}\Big\| D\k0(s)\Big\||\gaep_1(\omh,\kh)|,\nonumber
\end{align}
where
\begin{align*}
\left\|\Phi(\th;s)\right\|&\coloneq \max_{\th\in [\th_0,\th_0+1]}\max_{s\in [\th_0,\th]}\left|\Phi(\th;s)\right|,\\
\left\|\frac{d}{ds}\Phi(\th;s)\right\|&\coloneq \max_{\th\in [\th_0,\th_0+1]}\max_{s\in [\th_0,\th]}\left|\frac{d}{ds}\Phi(\th;s)\right|,
\end{align*}
and
\begin{equation}\label{M}
M\coloneq\frac{\|[Id-\Phi(\th_0+1;\th_0)]^{-1}\|\|\pp_{\th_0}\|}{\om_0}.
\end{equation}
We have used $[Id-\Phi(\th_0+1;\th_0)]^{-1}$ to denote the inverse of
$[Id-\Phi(\th_0+1;\th_0)]$ in the $(n-1)$-dimensional space $E$
introduced in Section \ref{sec:defspaces}.

Note that for the right hand sides of the inequalities \eqref{estga1}
and \eqref{estga2} above, each term is either quadratic in $a$,
$\beta_0$ or has a factor $\ep$. Under smallness assumptions of $a$,
$\beta_0$, and $\ep$, we will have $|\gaep_1(\omh,\kh)|\leq a$ and
$\|\gaep_2(\omh,\kh)\|\leq \beta_0$.

Now we consider the derivatives of $\gaep_2(\omh,\kh)$.

The first derivative 
  $\frac{d}{d\th}\Gamma^\ep_2(\omh,\kh)(\th)$ has the expression:
\begin{align*}
\Big(\frac{d}{d\th}&\Phi(\th;\th_0)\Big)\Big\{u_0+\frac{1}{\om_0}\int^{\th}_{\th_0}\Phi(s;\th_0)^{-1}\big(B^\ep(\kh,\omh,\ga,s)-\Gamma^\ep_1(\omh,\kh) D\k0(s)\big)ds\Big\}\\
&+\frac{1}{\om_0}\Big\{B^\ep(\kh,\omh,\ga,\th)-\Gamma^\ep_1(\omh,\kh) D\k0(\th)\Big\}.
\end{align*}
Recall that $\Phi(\th;\th_0)$ solves equation \eqref{var}. Therefore,
\begin{align*}
\Big\|\frac{d}{d\th}\gaep_2(\omh,\kh)\Big\|\leq& \frac{1}{\om_0}\|Df(\k0(\th))\|\|\gaep_2(\omh,\kh)\|+\frac{|\gaep_1(\omh,\kh)|}\Lip(Df) {\om_0}\Big\|D\k0(\th)\Big\|\\
&+\frac{1}{\om_0}\Big(\frac{1}{2}\Lip(Df) \|\kh\|^2+\ep \|\P(K,\om,\ga,\th)\|+|\omh|\Big\| D\kh(\th)\Big\|\Big)\\
\leq&\frac{1}{\om_0}\|Df(\k0(\th))\|\beta_0+\frac{a}{\om_0}\Big\|D\k0(\th)\Big\|\\
&\phantom{AAAA}+\frac{1}{\om_0}\Big(\frac{1}{2}\Lip(Df) \beta_0^2+\ep \|\P(K,\om,\ga,\th)\|+a\beta_1\Big).
\end{align*}
If $\ep$, $a$, and $\beta_0$ are small enough, we can choose $\beta_1$
to ensure that if $\| \frac{d}{d \th} \kh \| < \beta_1$, then
$\Big\|\frac{d}{d\th}\gaep_2(\omh,\kh)\Big\|<\beta_1$.

Now we proceed inductively, for $n\ge 2$,
$\frac{d^n}{d\th^n}\gaep_2(\omh,\kh)$ is an expression involving
$\Phi$, $\k0$, and their derivatives up to order $n$, as well as
$B^\ep$ and its derivatives up to order $n-1$. Within this expression,
$\k0$ and its derivatives are always multiplied by the small factor
$\Gamma^\ep_1(\omh,\kh)$, which has absolute value bounded by constant
$a$. It remains to consider $B^\ep$ and its derivatives.

Recall the definition of $B^\ep$ in \eqref{Bdef}, we now consider the
three terms in $B^\ep$ separately:
\begin{itemize}
\item 
For derivatives of $N(\th,\kh)$, we use the Faa di Bruno formula. The
$j$-th derivative of $N$ is an expression which contains derivatives
of $f$ up to order $j+1$, derivatives of $\kh$ up to order $j$.  All
the terms in derivatives of $N$ can be controlled taking advantage of
the fact that $N$ is of order at least $2$ in $\kh$.
\item
Derivatives of $\P$ are bounded thanks to the assumption
\eqref{H2.1}. Moreover, note that in $B^\ep$, $\P$ has the
perturbation parameter $\ep$ as its coefficient. Hence, this term is
less crucial.
\item
For the last term, $\omh D\kh(\th)$, its $j$-th derivative is $\omh
D^{j+1}\kh(\th)$. All are under control since $|\omh|<a$. Notice that
the $(n-1)$-th derivative of this term is $\omh D^{n}\kh(\th)$, and
this is the only place that $D^{n}\kh(\th)$ appears.
\end{itemize}

Taking all the terms above into consideration and using the triangle
inequality, we obtain bounds
\begin{equation}
\label{polyn}
\left\|\frac{d^n}{d\th^n}\gaep_2(\omh,\kh) \right\| \le P_n(a, \beta_0,\dotsc,\beta_{n-1})+\alpha\beta_n,
\end{equation}
where for each $n$, $P_n$ is a polynomial expression with positive
coefficients, and $\alpha<1$. The coefficients of $P_n$ are
combinatorial numbers multiplied by derivatives of $\k0$, $\P$, $f$,
and $\Phi(\th;\th_0)$. Therefore, we can choose recursively the
$\beta_i$'s such that right hand side of inequality \eqref{polyn} is
bounded by $\beta_n$.

Similar estimation can be obtained for the Lipschitz constant of
$\frac{d^l}{d\th^l}\gaep_2(\omh,\kh)$. Hence, we can choose $a$,
$\beta$ such that $\gaep\colon I_a\times\mathcal{B}_\beta \to
I_a\times\mathcal{B}_\beta $.
\end{proof}

\begin{remark}
Note that for $\ep$ small, we can get constant $a$ and each component
of $\beta$ are small.
\end{remark}

\begin{remark}\label{Schau}
Note that $I_a\times\mathcal{B}_\beta\subset\R\times C(\T,\R^n)$ is
compact and convex, and it is obvious that $\gaep\colon
I_a\times\mathcal{B}_\beta \to I_a\times\mathcal{B}_\beta $ is
continuous, so one could apply Schauder's fixed point Theorem to
obtain existence of the fixed point. Indeed, weaker assumptions than
assumption \eqref{H3.1} on $\P$ could also suffice to ensure
continuity of $\gaep$.

We will later prove that $\gaep$ is a contraction in $C^0$, which will
give local uniqueness of the fixed point and a-posteriori estimates on
the difference between an initial guess and the fixed point.

In principle, the Banach contraction theorem provides estimates of the
difference in $C^0$ norm, but, taking into account the propagated
bounds, we can use interpolation inequalities (Lemma~\ref{lem:inter} )
to obtain estimates in norms with higher regularity.  See
Section~\ref{sec:conc}.
\end{remark}

\subsection{Contraction Properties of $\gaep$}\label{sec:cont}
Define $C^0$-type distance on $I _a\times\mathcal{B} _\beta$:
\begin{equation}\label{dist}
d\left((\omh,\kh),(\omh',\kh')\right) \coloneq
\max\left\{|\omh-\omh'|,\|\kh-\kh'\|\right\}.
\end{equation}
\begin{lemma}\label{lem:contract}
For small enough $\ep$, $a$, and $\beta_0$(as in $\beta$), the
operator in \eqref{op} is a contraction on $I_a\times
\mathcal{B}_\beta$ with distance \eqref{dist}, i.e., there exists
$0<\mu<1$, such that
\begin{equation}\label{contract}
d\left(\gaep(\omh,\kh),\gaep(\omh',\kh')\right)<\mu\cdot
d\left((\omh,\kh),(\omh',\kh')\right).
\end{equation}
\end{lemma}
\begin{proof}
The proof of this lemma consists basically in adding and subtracting
and estimating by the mean value theorem.

We first list some useful inequalities for proving this lemma:
\[
\|N(\th,\kh)-N(\th,\kh')\|\leq \frac{1}{2}\Lip(Df) \big(\|\kh\|+\|\kh'\|\big)\|\kh-\kh'\|,
\]
where the norm $\|D^2 f\|$ is still interpreted as the norm in a
neighborhood of the periodic orbit of the unperturbed equation, as in
the proof of Lemma \ref{lem:prop}.

For $\th\in [\th_0,\th_0+1]$,
\begin{align*}
\left| \int^\th_{\th_0}\right.&\Phi(\th_0+1;s)\omh D\kh(s)ds-\left.\int^\th_{\th_0}\Phi(\th_0+1;s)\omh' D\kh'(s)ds\right|\\
&\leq \left(2\|\Phi(\th_0+1;\th)\|+\Big\|\frac{d}{d\th}\Phi(\th_0+1;\th)\Big\|\right)\left(\|\kh\||\omh-\omh'|+|\omh'|\|\kh-\kh'\|\right).
\end{align*}
Define 
\[
\om'=\om_0+\omh',\quad K'=\k0+\kh',
\]
similar to \eqref{kom}. 

By assumption \eqref{H3.1}, 
\begin{equation*}
\phantom{AAA}|\P(K, \om, \ga, \th)-\P(K', \om', \ga, \th)|\leq \alpha_{\rho,\delta}\max\left\{|\om-\om'|,\|K-K'\|\right\}.
\end{equation*}

Then,
\begin{align}\label{contga1}
|\gaep_1&(\omh,\kh)-\gaep_1(\omh',\kh')|\\
\leq& \frac{\|\pt_{\th_0}\|\Big(2\|\Phi(\th_0+1;\th)\|+\big\|\frac{d}{d\th}\Phi(\th_0+1;\th)\big\|\Big)}{\big| D\k0(\th_0)\big|}(\beta_0|\omh-\omh'|+a\|\kh-\kh'\|)\nonumber\\
&+\frac{\|\pt_{\th_0}\|\|\Phi(\th_0+1;\th)\|}{\big| D\k0(\th_0)\big|}\Big[\beta_0\Lip(Df) \|\kh-\kh'\|+\ep \alpha_{\rho,\delta} d\big((\omh,\kh),(\omh',\kh')\big)\Big].\nonumber
\end{align}

The initial conditions in both cases are:
\begin{align*}
u_0&=\frac{1}{\om_0}[Id-\Phi(\th_0+1;\th_0)]^{-1}\int^{\th_0+1}_{\th_0}\pp_{\th_0}\Phi(\th_0+1;s)B^\ep(\kh,\omh,\ga,s)ds;\\
u_0'&=\frac{1}{\om_0}[Id-\Phi(\th_0+1;\th_0)]^{-1}\int^{\th_0+1}_{\th_0}\pp_{\th_0}\Phi(\th_0+1;s)B^\ep(\kh',\omh',\ga,s)ds.
\end{align*}
As before, $[Id-\Phi(\th_0+1;\th_0)]^{-1}$ denotes the inverse of
$[Id-\Phi(\th_0+1;\th_0)]$ in the $(n-1)$-dimensional space $E$
introduced in Section \ref{sec:inveq}.

Therefore,
\begin{align}\label{contu0}
|u_0-&u_0'|\\
\leq& M\Big(2\|\Phi(\th_0+1;\th)\|+\Big\|\frac{d}{d\th}\Phi(\th_0+1;\th)\Big\|\Big)\big(\beta_0|\omh-\omh'|+a\|\kh-\kh'\|\big)\nonumber\\ 
&+M\|\Phi(\th_0+1;\th)\|\Big[\beta_0\Lip(Df) \|\kh-\kh'\|+\ep \alpha_{\rho,\delta}d\big((\omh,\kh),(\omh',\kh')\big)\Big],\nonumber
\end{align}
where $M$ is defined as in \eqref{M}. Therefore,
\begin{align}\label{contga2}
\|\gaep_2(\omh,\kh)-&\gaep_2(\omh',\kh')\|\\
\leq& \|\Phi(\th;\th_0)\||u_0-u_0'|+\frac{\|\Phi(\th;s)\|}{\om_0}\|D\k0(\th)\||\gaep_1(\omh,\kh)-\gaep_1(\omh',\kh')|\nonumber\\
&+\frac{\|\Phi(\th;\th_0)\|+\big\|\frac{d}{ds}\Phi(\th;s)\big\|+1}{\om_0}\big(\beta_0|\omh-\omh'|+a\|\kh-\kh'\|\big)\nonumber\\
&+\frac{\|\Phi(\th;s)\|}{\om_0}\Big[\beta_0\Lip(Df) \|\kh-\kh'\|+\ep \alpha_{\rho,\delta}d\big((\omh,\kh),(\omh',\kh')\big)\Big].\nonumber
\end{align}
Combining \eqref{contga1}, \eqref{contu0}, and \eqref{contga1}, if
$\ep$ is sufficiently small, $a$ and $\beta_0$ are chosen to be
sufficiently small, we can find $\mu$ such that \eqref{contract} is
true, $\gaep$ is a contraction.
\end{proof}

\subsection{Conclusion of the Proofs of Theorem \ref{thm:per}}\label{sec:conc}
There exists a fixed point $(\omh^*,\kh^*)$ of contraction
$\gaep$. According to Arzela-Ascoli Theorem (see Lemma \ref{closure}
in Appendix), $(\omh^*,\kh^*)\in I_a\times\mathcal{B}_\beta $, hence
is a solution of the functional equation \eqref{inveq} with desired
regularity. Then, $K=\k0+\kh^*$ gives a parameterization of the
periodic orbit of \eqref{sdde}.

The proof based on fixed point approach leads to a-posteriori type of
results. Suppose we start with initial guess $(\omh^0,\kh^0)$ for
$(\omh,\kh)$, since $\gaep$ is contractive, see equation
\eqref{contract}, we have
\begin{equation}
d\left((\omh^0,\kh^0),(\omh^*,\kh^*)\right)<\frac{1}{1-\mu} d\left((\omh^0,\kh^0),\gaep(\omh^0,\kh^0)\right).
\end{equation}
Therefore, if we have a good choice of initial guess such that the
error in the fixed point equation,
$d\big((\omh^0,\kh^0),\gaep(\omh^0,\kh^0)\big)$, is small, then we
know the fixed point is close to the initial guess.

Using interpolation inequalities in Lemma \ref{lem:inter}, we also have
\[
\|\kh^0-\kh^*\|_{C^m} \leq C\|\kh^0-\kh^*\|^{1-\frac{m}{\ell}}_{C^0}
\]
for $0\leq m<\ell$, where the constant $C$ depends on $\ep,~m,~\ell$,
and $\beta_i,~ i=0,1,\dotsc,\ell$. In particular, the distance between
the initial guess $(\omh^0,\kh^0)=(0,0)$ and the fixed point
$(\omh^*,\kh^*)$ is of order $\ep$, therefore, $\|\kh^*\|_{C^m}$ is
small for $0\leq m<\ell$.  This finishes the proof of Theorem
\ref{thm:per}.

\subsection{Comments on proof of Theorem \ref{thm:smooth}}\label{sec:cosm}

A very similar method proves Theorem \ref{thm:smooth}. Now we view
$\omh$ as a function of $\ga$, and $\kh$ as a function of $\th$ and
$\ga$. Define operator $\widetilde\Gamma^\ep$ of the same format as in
\eqref{op1} and \eqref{op2} on the space
$\pazocal{I}\times\pazocal{F}$, where $\pazocal{I}$ contains
$C^{\ell+\Lip}$ functions from set $O$ to $I$ and $\pazocal{F}$
contains $C^{\ell+\Lip}$ functions from $\T\times O$ to $\R^n$, with
bounded derivatives similar to \eqref{space}. We can then prove that
for small enough $\ep$, and suitable choices for bounds of
derivatives, $\widetilde\Gamma^\ep$ maps
$\pazocal{I}\times\pazocal{F}$ to itself using assumption
\eqref{H2.2}, and is a contraction in $C^0$ norm, taking advantage of
assumption \eqref{H3.2}. Therefore, there exists a fixed point for
$\widetilde\Gamma^\ep$ in the space $\pazocal{I}\times\pazocal{F}$
solving equation \eqref{inveq}. Same as above, Theorem
\ref{thm:smooth} is proved.

\section{Delay Perturbation to Autonomous ODE}\label{sec:examauto}

In this section we show how several concrete examples fit into our
general result Theorem~\ref{thm:per}.  In all the cases, we will show
how to construct the operators $\P$ and to verify the properties in
assumptions (H2) and (H3).

\subsection{State-Dependent Delay Perturbation}
\label{sec:sdde}

An important class of equations that one can consider is DDEs with
state-dependent delays (backward or forward or mixed):
\begin{equation}\label{sdde}
\dot x(t)=f(x(t))+\ep P\Big(x(t), x\big(t-r(x(t))\big),\gamma\Big),
\end{equation}
where $P\colon \R^n\times\R^n\times O\to \R^n$ is a $C^{\infty}$ map,
$r\colon \R^n\to [-h,h]$ is $C^{\infty}$, $h$ is a positive constant.

Note that in this case, the operator $\P$ is,
\begin{equation}\label{sddep}
\P(K, \om, \ga, \th)=P(K(\th), \kt(\th), \ga)
\end{equation}
where $\kt(\th)\coloneq K(\th-\om r(K(\th)))$ is caused by the delay.

\begin{remark} 
Note that the operator $\P$ involves the composition operator, whose
differentiability properties are very complicated (See \cite{LO99} for
a systematic study). Hence, using the standard strategy of studying
variational equations etc. to study regularity of the evolution will
be rather complicated. Indeed, it will be hard to go beyond the first
derivative.

On the other hand, the present strategy, only requires much simpler
results.  We only need to get bounds on the derivatives of $\P$
assuming bounds on the derivatives of $K$.
\end{remark} 

Applying the composition Lemma \ref{lem:com} repeatedly, we know that
the $\P$ above satisfies \eqref{H2.1}. With the standard adding and
subtracting terms method, one gets that $\P$ satisfies
\eqref{H3.1}. Similarly, $\P$ satisfies \eqref{H2.2} and
\eqref{H3.2}. Thus, Theorem \ref{thm:per} and \ref{thm:smooth} can be
applied.
 
Note also that for the above equation \eqref{sdde}, we are able to
prove that the operator $\gaep$ is a contraction under
$C^{\ell-1+\Lip}$ norm in the second component, by using
Lemma~\ref{compositionLipschitz}.

We can improve the regularity conclusion of Theorem \ref{thm:per} for
this case. Indeed, thanks to the high regularity of $P$ and $r$ in
equation \eqref{sdde}, once we have that the parameterization $K$ of
the periodic orbit is $C^1$ in $\th$, we can use the standard
bootstrapping argument to conclude that $K$ is $C^\infty$.

We can also consider more general state-dependent delays:
\begin{equation}\label{sdde2}
\dot x(t)=f(x(t))+\ep P\Big(x(t),
x\big(t-r\left(x_t,\ga\right)\big),\gamma\Big),
\end{equation}
where $r\colon \mathcal{R}[-h,h]\times O \to \R$, positive constant
$h$ is an upper bound for $|r|$.

In this case,
\begin{equation}
\P(K, \om, \ga, \th)=P\left(K(\th), K\big(\th-\om
r(K_{\th,\om},\ga)\big), \ga\right),
\end{equation}
where $K_{\th,\om}\colon[-h,h]\to \R^n$ is defined by
\begin{equation}\label{kthetaom}
K_{\th,\om}(s)\coloneq K(\th+\om s).
\end{equation}

If $r$ is chosen such that (H2) and (H3) are verified, Theorems
\ref{thm:per} and \ref{thm:smooth} can be applied.

\subsection{Distributed Delay Perturbation}

Our results apply to models with distributed delays as well
\begin{equation}\label{distribu}
\dot x(t)=f(x(t))+\ep P\Big(x(t), \int^0_{-r} x_t(s)d\mu(s),\gamma\Big),
\end{equation}
where $P\colon \R^n\times\R^n\times O\to \R^n$ is a $C^{\infty}$ map,
$r$ is a constant, and $\mu$ is a signed Borel measure.  In this case,
\begin{equation}
\P(K, \om, \ga, \th)=P\left(K(\th), \int^0_{-r} K_{\th,\om}(s)d\mu(s),
\ga\right),
\end{equation}
where $K_{\th,\om}$ is defined in \eqref{kthetaom}.

Above $\P$ verifies (H2), since we only care about derivatives with
respect to $\th$. It is not hard to see that (H3) is also satisfied in
this case. Therefore, Theorems \ref{thm:per} and \ref{thm:smooth}
apply.

\subsection{Remarks on further applicability of Theorem~\ref{thm:per}}

\begin{remark} 
It is straightforward to see that our results could be applied to
systems similar to above systems with multiple forward or backward
delays.
\end{remark}

\begin{remark}
In some applications, the delays are defined by some implicit
relations from the full trajectory.

The Theorem~\ref{thm:per} can be applied if we can justify
\eqref{H2.1}, \eqref{H3.1}, \eqref{H2.2}, and \eqref{H3.2}.  Notice
that we only need to justify these hypothesis when the $\kh$ lies in
ball in a space of differentiable functions. In such a case, we can
often just use the implicit function theorem.
\end{remark}

\begin{remark} 
The results so far do not include the models in which the perturbation
is just adding a small delay.  This small delay perturbation is more
singular and seems to require extra assumptions and slightly different
proofs.  The extension of the results to the small delay case is done
in Section~\ref{sec:small_delays}.
\end{remark}

\begin{remark} \label{rmk:bootana}
In the case of state dependent delay or distributed delay with a
smooth $f$ and $r$, it is automatic to show that if $K$ is $C^\ell$,
the right hand side of \eqref{eqk} is $C^\ell$, hence, looking at the
left hand side of \eqref{eqk}, $K$ is $C^{\ell +1}$. The bootstrap
stops only when we do not have any more regularity of $f$ or $r$.

So, in case that $f$ and $r$ are $C^\infty$, we obtain that the $K$
is $C^\infty$.

One natural question that deserves more study is whether in the case
that $f$ and $r$ are analytic, the $K$ is analytic.  The remarkable
paper \cite{MNAna} contains obstructions that show that equations with
time dependent delays -- heuristically better behaved than the ones
considered here, may fail to have analytic solutions. In view of these
results, it is natural to conjecture that the periodic solutions
produced here, could fail to be analytic even if $f$ and $r$ are
analytic.
\end{remark}

\section{Delay Perturbation to Non-autonomous Periodic ODE}\label{sec:examper}

Time periodic systems appear in many problems in physics, for example,
see Section \ref{sec:electrodynamics}. And when there are conserved
quantities in the ODE systems, periodic orbits cannot satisfy the
assumption \eqref{H1}. These are the motivations to consider a
non-autonomous ODE:
\begin{equation}\label{odenon}
\dot x(t)=f(x(t),t),
\end{equation}
where $f\colon\R^n\times \frac{1}{\om_0} \T\to \R^n$ ($f$ is periodic
in $t$ with period $\frac{1}{\om_0}$). Add the perturbation:
\begin{equation}\label{fdenon}
\dot x(t)=f(x(t),t)+\ep P(x_t,\ga),
\end{equation}

Using the standard method of adding an extra variable to equation
\eqref{odenon} to make it autonomous, we will see we can reduce the
problem to the previous case. The autonomous equation corresponding to
\eqref{odenon} is
\begin{equation}\label{odeauto}
\begin{pmatrix}
\dot x(t)\\
\dot \tau(t)
\end{pmatrix}=g(x,\tau)\coloneq \begin{pmatrix}
f(x,\tau)\\1
\end{pmatrix}.
\end{equation}

Denote $\Psi$ as the solution of the variational equation for the
periodic orbit of \eqref{odeauto}:
\begin{equation}
\om_0\frac{d}{d\theta}\Psi(\theta;\theta_0)=Dg\left(\k0(\theta),\frac{\th}{\om_0}\right)\Psi(\theta;\theta_0),\quad \Psi(\theta_0;\theta_0)=Id.
\end{equation}
Since
\begin{equation*}
Dg\left(\k0(\theta),\frac{\th}{\om_0}\right)=\begin{pmatrix}
D_1f\left(\k0(\th),\frac{\th}{\om_0}\right)&D_2f\left(\k0(\th),\frac{\th}{\om_0}\right)\\
0&0
\end{pmatrix},
\end{equation*}
we have
\begin{equation}
\Psi(\theta;\theta_0)=\begin{pmatrix}
\Phi(\theta;\theta_0)&*\\
0&1
\end{pmatrix},
\end{equation}
where 
\begin{equation}\label{periodicvar}
\om_0\frac{d}{d\theta}\Phi(\theta;\theta_0)=D_1f\left(\k0(\th),\frac{\th}{\om_0}\right)\Phi(\theta;\theta_0).
\end{equation}
If $\Psi$ satisfies assumption \eqref{H1}, then $\Phi(\th_0+1;\th_0)$
has no eigenvalue 1.

Equivalently, we could start the discussion in this section directly
with the following assumption on $\Phi$ defined in
\eqref{periodicvar}:
\begin{enumerate}
\renewcommand{\theenumi}{H\arabic{enumi}''}
\renewcommand{\labelenumi}{(\theenumi)}
 \item \label{H1''} $\Phi(\th_0+1;\th_0)$ has no eigenvalue 1.
\end{enumerate}

Under either assumption \eqref{H1} on $\Psi$ or assumption
\eqref{H1''} on $\Phi$, we are able to solve the invariance equation
\eqref{eqk} without adjusting the frequency. More precisely,
\eqref{eqk} becomes:
\begin{equation}
\om_0DK(\th)=f\left(K(\th),\frac{\th}{\om_0}\right)+\ep \P(K,\om_0,\ga,\th),
\end{equation}
Let $K=\k0+\kh$ as in \eqref{kom}, we are led to 
\begin{equation}
\om_0D\kh(\th)-D_1f\left(\k0(\th),\frac{\th}{\om_0}\right)\kh(\th)=B^\ep(\th,\om_0,\kh,\gamma),
\end{equation}
where
\begin{align*}
&B^\ep(\th,\om_0,\kh,\gamma)\coloneq N(\th, \kh)+\ep \P(K,\om_0,\ga,\th),\\
&N(\th, \kh)\coloneq f\left(K(\th),\frac{\th}{\om_0}\right)-f\left(\k0(\th),\frac{\th}{\om_0}\right)-D_1f\left(\k0(\th),\frac{\th}{\om_0}\right)\kh(\th).
\end{align*}

Now we can define an operator $\Upsilon^\ep$ on the space $\B _\beta$
(see \eqref{space}) very similar to the second component of $\gaep$
introduced in section \ref{sec:operator}.
\begin{equation}\label{periodicop}
\Upsilon^\ep(\kh)(\th)\coloneq\Phi(\th;\th_0)u_0+\frac{1}{\om_0}\int^{\th}_{\th_0}\Phi(\th;s)B^\ep(s,\om_0,\kh,\gamma)ds,
\end{equation}
where
\begin{equation}
u_0=\frac{1}{\om_0}[Id-\Phi(\th_0+1;\th_0)]^{-1}\int^{\th_0+1}_{\th_0}\Phi(\th_0+1;s)B^\ep(s,\om_0,\kh,\gamma)ds.
\end{equation}
We have employed that, in the periodic case, the matrix
$[Id-\Phi(\th_0+1;\th_0)]$ is invertible.

Under the assumption that $\P$ satisfies \eqref{H2.1}, \eqref{H3.1},
\eqref{H2.2} and \eqref{H3.2}, we can prove that $\Upsilon^\ep$ has a
fixed point $\kh^*\in\B _\beta$ by proving $\Upsilon^\ep\colon\B
_\beta\to\B _\beta$ (similar to Lemma \ref{lem:prop}) and
$\Upsilon^\ep$ is a contraction (similar to Lemma
\ref{lem:contract}). The periodic orbit of \eqref{fdenon} is
parameterized by $K=\k0+\kh^*$.  The analysis of the operator
$\Upsilon^\ep$ in \eqref{periodicop} is actually simpler that the
analysis presented for the operator $\Gamma^\ep$ in \eqref{op} because
we do not need to adjust the frequency.

\begin{remark}
Similarly, we can also consider a non-autonomous perturbation $P(t,
x_t,\ga)$, we need that $P$ to be periodic in $t$ with the same period
$\frac{1}{\om_0}$.
\end{remark}

\section{The case of small delays} 
\label{sec:small_delays}
Many problems in the literature lead to 
equations of the form:
\begin{equation} \label{smalldelay}
  \begin{split} 
    &y'(t) = g(y(t -\ep r) )\\
    &y'(t) = f(y(t -\ep r),t)
  \end{split}
\end{equation}
where $r$ could be either a constant, an explicit function of $t$, a
function of $y(t)$, or $y_t$, or an implicit function, and may depend
on $\ep$; and $f$ is periodic of period $1$ in $t$. Indeed, our
results apply also to variants of \eqref{smalldelay} with
perturbations involving several forward or backward delays.

In problems which present feedback loops, the feedback takes some time
to start acting.  The problems \eqref{smalldelay} correspond to the
feedback taking a short time to start acting.

Equations of the form \eqref{smalldelay} play an important role in
electrodynamics, where the small parameter $\ep = \frac{1}{c}$ is the
inverse of the speed of light and the delay $r$ is a functional that
depends on the trajectory.  Given the physical importance of
electrodynamics, we devote Section~\ref{sec:electrodynamics} to give
more details and to show that it can be reduced to
Theorem~\ref{thm:smalldelay}.

Introducing a small delay to the ODE is a very singular perturbation,
since the phase space becomes infinite dimensional. The limit is
mathematically harder because the effect of a small delay is similar
to adding an extra term containing the derivative $y(t -\ep r) \approx
y(t) +\ep y'(t)r $.  This shows that, heuristically, the perturbation
is of the same order as the equation.

\begin{remark} 
In the physical literature, one can find the use of higher order
expansions to obtain heuristically even higher order equations, see
\cite{Dirac38}. As a general theory for all the solutions of the
equations, these theories have severe paradoxes
(e.g. \emph{preacceleration}). The results of this paper show, however
that the non-degenerate periodic solutions produced in many of these
expansions, since they are very approximate solutions of the
invariance equation, approximate true periodic solutions of the full
system.
\end{remark} 

As a reflection of the extra difficulty of the small delay problem
compared with the previous ones, the main result of this section,
Theorem~\ref{thm:smalldelay}, requires a more delicate proof than
Theorem~\ref{thm:per} and we need stronger regularity to obtain the
$C^0$ contraction.

An important mathematical paper on the singular problem of small delay
is \cite{Chiconedelay}.  We also point out that, there is a
considerable literature in the formal study of $\frac{1}{c}$ limit in
electrodynamics and in gravity \cite{LandauL,
  Plass61,Bel71,Rohrlich81}.  Many famous consequences of relativity
theory (e.g. the precession of the perihelion of Mercury) are only
studied by formal perturbations.

Formal expansions of periodic and quasiperiodic solutions for small
delays were considered in \cite{CasalCL20}.  The results of this
section establish that the formal expansions of periodic orbits
obtained in \cite{CasalCL20} correspond to true periodic orbits and
are asymptotic to the true periodic solutions in a very strong sense.

In this section, we establish results on persistence of periodic
orbits for the models in \eqref{smalldelay}, see
Theorem~\ref{thm:smalldelay}.  As we will see, when we perform the
detailed discussion, we will not be able to reduce
Theorem~\ref{thm:smalldelay} to be a particular case of
Theorem~\ref{thm:per}.  The proof of Theorem~\ref{thm:smalldelay} will
be very similar to that of Theorem~\ref{thm:per} and which is based on
the study of operator $\Gamma^\ep$ very similar to those in
\eqref{op}.  Nevertheless, the analysis of the operator $\Gamma^\ep$
in the current case will require to take advantage of an extra
cancellation.

\subsection{Formulation of the results} 

Our main result for the small delay problem \eqref{smalldelay} is as
follows. Without specifying the delay functional $r$, we will use
$r(\om,K,\ep)$ to denote the expression after substituting $K(\th+\om
t)$ into $r$ and letting $t=0$.

\begin{theorem} \label{thm:smalldelay} 
For integer $\ell\ge 3$, assume that the function $g$ (resp. $f$) in
\eqref{smalldelay} is $C^{\ell + \Lip}$.

Assume that for $\ep = 0 $, the ordinary differential equation $y' =
g(y) $ has a periodic orbit satisfying \eqref{H1}.(resp.  $y' = f(y,t)
$ has a periodic orbit satisfying \eqref{H1''} ). We denote by $K_0$
the parameterization of this periodic orbit with frequency $\om_0$.

Recall that $U_\rho$ is the ball of radius $\rho$ in $C^{\ell +
  \Lip}(\T,\R^n)$ centered at $K_0$, and $B_\delta$ is the interval
with radius $\delta$ centered at $\om_0$.

Recall distance $d$ defined in \eqref{dist}. Assume that the delay
functional $r$ satisfies
\begin{equation}\label{rbounds} 
\|r(\om,K,\ep) \|_{C^{\ell-1+\Lip}(\T,\R^n)} \le
\phi_{\rho,\delta}(\ep) \quad \forall K \in U_\rho,~ \om\in B_\delta
\end{equation}
for some $\phi_{\rho,\delta}(\ep)>0$, with
$\ep\phi_{\rho,\delta}(\ep)\to 0$ as $\ep\to 0$.  And for any $K_1,
K_2 \in U_\rho$, $ \om_1, \om_2\in B_\delta$, there is
$\alpha_{\rho,\delta}(\ep)>0$, with $\ep\alpha_{\rho,\delta}(\ep)\to
0$ as $\ep\to 0$, such that
\begin{equation}\label{c0for_r}
  \| r(\om_1, K_1,\ep) - r(\om_2, K_2,\ep) \|_{C^0} \le
  \alpha_{\rho,\delta}(\ep) d\big((\om_1,K_1),(\om_2,K_2)\big).
\end{equation}

Then, there exists $\ep_0 > 0$ such that for $\ep \le \ep_0$, the
problem \eqref{smalldelay} admits a periodic solution. There is a
parameterization $K$ of the periodic orbit which is close to $K_0$ in
the sense of $C^{\ell}$.
\end{theorem} 

\begin{remark} \label{increaseregularity} 
As before, the requirements of smallness in $\ep$ for
Theorem~\ref{thm:smalldelay} depend on the regularity considered.

In many applied situations, the $g$ and $f$ considered are $C^\infty$
or even analytic.  (for example in the electrodynamics applications
considered in Section~\ref{sec:electrodynamics}).  In such a case, we
can consider any $\ell$ by assuming $\ep$ is small enough.

This allows us to obtain the a-posteriori estimates in more regular
spaces as $\ep$ goes to zero.

Hence, the formal power series in \cite{CasalCL20} are asymptotic in
the strong sense that the error in the truncation is bounded by a
power of $\ep$, where a stronger norm can be used for smaller $\ep$.
\end{remark} 

We leave for the reader the formulation of a corresponding result for
the smooth dependence on parameters similar to Theorem
\ref{thm:smooth}.  The proof requires only small modifications from
discussion in Section \ref{smalldelayproof}, see comments in Section
\ref{sec:cosm}.

The proof of Theorem~\ref{thm:smalldelay} will be given in Section
\ref{smalldelayproof}.  We first find the operator in this case. Then
for the operator, we prove Lemma \ref{lem:prop} in Section
\ref{sec:propsmalldelay}, and prove Lemma \ref{lem:contract} in
Section~\ref{contractionsmalldelay}. The existence of fixed point of
the operator is thus established. As it turns out, the analysis of the
operator requires more care than in the case of Theorem~\ref{thm:per}.

\subsection{Formulating the existence of a fixed point operator} \label{smalldelayproof}

The equations \eqref{smalldelay} can be rearranged as
\begin{equation}\label{smalldelayrearranged}
  \begin{split} 
 y'(t) &= g( y(t) ) +
\left[  g( y(t - \ep r) ) - g( y(t)) \right] \\
&= g( y(t) ) -   \ep\int_0^1 \left[ D g\big( y(t - s \ep r) \big) Dy( t -s \ep r) r\right]ds \\
 y'(t) &= f(t, y(t) ) +  
\left[  f(t, y(t - \ep r) ) - f(t, y(t)) \right] \\
&= f(t, y(t) ) -   \ep\int_0^1  \left[  D_1 f\big(y(t - s \ep r),t \big) Dy(t -s \ep r)  r\right] ds
\end{split}
\end{equation}

For typographical convenience, we will discuss only the autonomous
case, which is the most complicated. We refer the reader to Section
\ref{sec:examper} to see how the discussion simplifies in the periodic
case (the most relevant case for applications to electrodynamics).

Note that \eqref{smalldelayrearranged} is in the form of \eqref{fde},
with the operator $P$ defined as
\begin{equation} \label{Pdefined} 
P[y](t) \coloneq-\int_0^1 \left[ D g\big( y(t - s \ep r) \big) Dy( t -s \ep r) r\right]\, ds.
\end{equation}
Then,
\begin{equation}\label{CPdefined}
\P(K, \om, \ga, \th)\coloneq-\int_0^1 \left[ D g\big( K(\th - \ep s \om r) \big) DK(\th -\ep s \om r) \om r\right]ds,
\end{equation}
where the $r$'s are $r(\om, K)$, the delay functional evaluated on the
periodic orbit.

We define operator $\Gamma^\ep$ in the same way as in
Section~\ref{sec:proof}, substituting the expression of $\P$ in
\eqref{CPdefined} into the general formula in \eqref{op}.

In this section, we will proceed as before and show Lemmas
\ref{lem:prop} and \ref{lem:contract} are true for the resulting
operator $\gaep$ with $\P$ defined in \eqref{CPdefined}.

Lemma \ref{lem:prop} is proven in this case, same as above, by
noticing $\P$ satisfies assumption \eqref{H2.1}. The proof for Lemma
\ref{lem:contract} is slightly different from before. In Section
\ref{sec:proof}, we only needed to take advantage of the Lipschitz
property of the operator $\P$ (assumption \eqref{H3.1}). In the
present case, we will have to take into account that the operator
$\Gamma^\ep$ involves not only $\P$, but also an integral, which has
nice properties that compensate the bad properties of $\P$.

\subsubsection{Propagated bounds} \label{sec:propsmalldelay} 
We observe that if $ K \in U_\rho$ and $\om\in B_\delta $, by the
assumption \eqref{rbounds}, $r(\om, K,\ep)$ is in a $C^{\ell -1+
  \Lip}$ ball of size $\phi_{\rho,\delta}(\ep)$ and, using the
estimates on composition, Lemma~\ref{lem:com}, so is $K(t - \ep s \om
r )$.  If $g\in C^{\ell + \Lip}$, then $Dg \in C^{\ell -1 + \Lip}$ and
we conclude that $Dg\circ K(t - \ep s \om r )$ is contained in a
$C^{\ell -1 + \Lip}$ ball.

We also have that if $K$ is in a $C^{\ell + \Lip}$ ball, $DK( t - \ep
s \om r) \om r $ is in a $C^{\ell -1 + \Lip}$ ball whose size is a
function of $\rho$ and $\phi_{\rho,\delta}(\ep)$.

Putting it all together we obtain that \eqref{H2.1} is true for $\P$
defined in \eqref{CPdefined}. Therefore, Lemma \ref{lem:prop} is
proven in this case.

\subsubsection{Contraction in $C^0$} \label{contractionsmalldelay}
Before estimating $\gaep(\omh,\kh)-\gaep(\omh',\kh')$, we estimate
$\P(K, \om, \ga, \th)-\P(K', \om', \ga, \th)$ (we denote by $r, r'$
the two delay terms corresponding to $\om, K$, and $\om', K'$
respectively).

As usual, adding and subtracting, we obtain that the difference in the
integrands in $\P$,
\[
D g\big( K(\th - \ep s \om r) \big) DK(\th -\ep s \om r) \om r-D g\big( K'(\th - \ep s \om' r') \big) DK'(\th -\ep s \om' r') \om' r'
\]
can be written as a sum of 8 differences in which only one of the
objects changes, see \eqref{breakup} below.  As it turns out, 7 of
them will be straightforward to estimate and only one of them will
require some effort.  We give the details.

\begin{equation}  \label{breakup} 
\begin{split} 
&[D g\big( K(\th - \ep s \om r) \big)- D g\big( K'(\th - \ep s \om r) \big)]DK(\th -\ep s \om r) \om r \\
&+[D g\big( K'(\th - \ep s \om r) \big)- D g\big( K'(\th - \ep s \om' r) \big)] DK(\th -\ep s \om r) \om r \\ 
&+[D g\big( K'(\th - \ep s \om' r) \big)- D g\big( K'(\th - \ep s \om' r') \big)] DK(\th -\ep s \om r) \om r \\
&+D g\big( K'(\th - \ep s \om' r') \big) [DK-DK'](\th -\ep s \om r)\om r \\ 
&+D g\big( K'(\th - \ep s \om' r') \big) [DK'(\th -\ep s \om r)-DK'(\th -\ep s \om' r)]\om r\\ 
&+D g\big( K'(\th - \ep s \om' r') \big) [DK'(\th -\ep s \om' r)-DK'(\th -\ep s \om' r')]\om r\\ 
&+D g\big( K'(\th - \ep s \om' r') \big) DK'(\th -\ep s \om' r') (\om-\om' ) r \\ 
&+D g\big( K'(\th - \ep s \om' r') \big) DK'(\th -\ep s \om' r') \om' (r-r') 
\end{split} 
\end{equation} 
 
All the terms except for the 4th term are straightforward to estimate
in $C^0$ by some constant multiple of
$d\big((\omh,\kh),(\omh',\kh')\big)$, keeping in mind bounds on the
$C^{\ell+\Lip}$ norms of $g$, $K$, $K'$, and $r$ \big(see assumption
\eqref{rbounds}\big), and the assumption~\eqref{c0for_r}. We consider
the first term for an example, the rest is similar.

\begin{equation}
\begin{split} \label{breakupestimate}
\|[D g\big( K(\th - \ep s \om r) \big)- D g\big( K'(\th - \ep& s \om r) \big)]DK(\th -\ep s \om r) \om r\|_{C^0}\\
&\le \om\|D^2g\|\|DK\|\|r\|_{C^0} \|\kh-\kh'\|_{C^0}
\end{split}
\end{equation}

Observe the form of the operator $\gaep$ in \eqref{op}. Note that if
we have a bound of
\[
\int^{\th}_{\th_0}\Phi(\th;s)\big(\P(K, \om, \ga, s)-\P(K', \om', \ga, s)\big)ds.
\] 
by a multiple of $d\big((\omh,\kh),(\omh',\kh')\big)$, we prove
Lemma~\ref{lem:contract}.
\medskip

All terms except the 4th one in \eqref{breakup} are controlled using
estimates similar to \eqref{breakupestimate}. Hence, to complete the
proof, we just need to estimate the part coming from the 4th term in
\eqref{breakup}.  We will take advantage of the integral which is an
operator that improves the bounds.

We use integration by parts to get: 
\[
\begin{split} 
&\int^{\th}_{\th_0} \Phi(\th;s )\int^1_0 D g\big( K'(s - \ep \tau \om' r') \big) [DK-DK'](s -\ep \tau \om r)\om rd\tau ds  \\
&= \int_0^1 \int^{\th}_{\th_0} \Phi(\th;s)Dg(\cdot)\frac{\om r}{1 - \ep \tau \om \frac{dr}{ds}}
(1 - \ep \tau \om \frac{dr}{ds})[DK-DK'](s -\ep \tau \om r) ds d\tau\\ 
&= \int_0^1
\Big[  \Phi(\th;s)Dg(\cdot)\frac{\om r}{1 - \ep \tau \om \frac{dr}{ds}} [K-K'](s -\ep \tau \om r)\big|_{s=\th_0}^{s = \th} \\ 
&\phantom{AAA} -  \int^{\th}_{\th_0} \frac{d}{ds} \Big(\Phi(\th;s)Dg(\cdot)\frac{\om r}{1 - \ep \tau \om \frac{dr}{ds}}\Big)  
[K-K'](s -\ep \tau \om r) ds
\Big]d\tau
\end{split} 
\]

The $C^0$ norm of the above expression is bounded by a multiple of $
\|K-K'\|_{C^0}\le d\big((\omh,\kh),(\omh',\kh')\big)$, we have proved
Lemma \ref{lem:contract} in this case. The proof of Theorem
\ref{thm:smalldelay} is finished.
\begin{remark}
Note that we need to differentiate $r$ along the periodic orbit twice
in the above expression, that is why we required $\ell\ge 3$ in
Theorem \ref{thm:smalldelay}, so that $r(\om,K,\ep)$ is more than
$C^2$.
\end{remark}

\section{Delays implicitly defined by the solution. Applications 
to electrodynamics}
\label{sec:electrodynamics} 

In this section, we show how to deal with delays that depend
implicitly on the solution. The main motivation is electrodynamics, so
we deal with this case in detail, but we formulate a more general
mathematical result in Section~\ref{sec:mathematicalresult}.

We point out that implicitly defined delays appear naturally in other
problems in which the delay of the effect is related to the state of
the system.  As we indicate later, the explicit state dependent delays
appeared in \eqref{sdde} are often approximations of implicitly
defined delays. One corollary of our treatment is a justification of
the fact that the periodic solutions of this approximation are an
approximation to the true periodic solutions.

\subsection{Motivation from Electrodynamics} 
One of the original motivations for the whole field of delay equations
was the study of forces in electrodynamics. The forces among charged
particles, depend on the positions of the particles.  Since the
signals from a particle take time to reach another particle, this
leads to a delay equation. Notice that the delay depends on the
position (at a previous time) so that the delay is obtained by an
implicit equation on the trajectory.  This formulation was proposed
very explicitly in \cite{WheelerF49}, which we will follow.

\begin{remark} 
An alternative description of electrodynamics uses the concept of
fields.  One problem of the concept of fields is to explain why
particles do not interact with their own fields. We refer to
\cite{Spohn04} for a very lucid physical discussion of the paradoxes
faced by a coherent formulation of classical electrodynamics.
\end{remark} 

\begin{remark} 
Many Physicists object to \cite{WheelerF49} that it does not make
clear what is the phase space and what are the initial conditions.

In this paper, we show that one does not need to answer these question
to construct a theory of periodic solutions.  We hope that similar
results hold for other types of solutions.  So that one can have a
systematic theory of many solutions that resemble the classical ones.

Of course, it should also be possible to construct other solutions
that are completely different from those of the systems without
delays.

\end{remark}

\begin{remark} 
Even if one can have a rich theory of perturbative solutions, It is
not clear that these solutions fit together in a smooth manifold. The
paper \cite{CasalCL20} develops asymptotic expansions, which suggests
that the resulting solutions may be difficult to fit together in a
manifold.

We speculate that this may give a way to reconcile the successes of
\emph{predictive mechanics} \cite{Bel71} with the no-interaction
theorems \cite{CurrieS}. It could well happen that the results of
predictive mechanics apply to the abundant solutions we construct,
but, according to the no-interaction theorem, this set cannot be all
the initial conditions.  Of course, these speculations are far from
being theorems.
\end{remark} 

\subsection{Mathematical formulation}

If we consider (time-dependent) external and magnetic fields as
prescribed, the equations of a system of $N$ particles in $\R^3$ are,
denoting by $q_i(t)$ the position of the $i$-th particle.
\begin{equation}\label{electrodynamics} 
  q''_i(t) = A_{\text{ext}}(t, q_i(t), q'_i(t)) + \sum_{j \ne i}
  A_{i,j}\left( q_i(t), q'_i(t), q_j(t - \tau_{ij}), q'_j(t -
  \tau_{ij} )\right)
\end{equation}
where the time delay is defined implicitly by ($c$ is the speed of
light)
\begin{equation}\label{delaydefined} 
  \tau_{ij}(t)  = \frac{1}{c} \big| q_i(t) -   q_j(t -\tau_{ij}(t) ) \big|.
\end{equation} 
For more explicit expressions, we refer to \cite{WheelerF49, Rohrlich,
  Driver}.  We just remark that \eqref{electrodynamics} is the usual
equation of acceleration equals force divided by the mass.  The
relativistic mass has some complicated expression depending on the
velocity.

The term $A_{\text{ext}}$ denotes the external force. The terms
$A_{ij}$ correspond to the Coulomb and Lorenz forces of the fields
obtained from Li\'enard–Wiechert potentials. This is a standard
calculation which is classical in electrodynamics, see \cite{LandauL,
  Jackson, Zangwill}. Roughly, they are the Coulomb and Ampere
(electric and magnetic) forces at previous times but some derivative
terms appear.

We observe that \eqref{electrodynamics} is in the form imposed by the
principle of relativity, and that any force which is relativistically
invariant should have the form \eqref{electrodynamics} with, of
course, different expressions for the terms $A_{ij}$. Hence, the
treatment discussed here should apply not only to electrodynamics but
also to any forces subject to the rules of special relativity.

The exact form of the equations does not play an important role in
this paper.  We point out some properties that play a role:

\begin{enumerate}
\item
  The expressions defining the forces are algebraic expressions. They
  have singularities when there are collisions ($q_i(t) = q_j(t)$ for
  some $i\ne j$) or when some particle reach the speed of light
  ($|q'_i(t) | = c$ for some $i$).
\item 
The delays $\tau_{ij}$ as in \eqref{delaydefined} are subtle. The
expression of $\tau_{ij}$ involve a small parameter $\ep\coloneq 1/c$,
and the delays can be approximated in first order as:
  \begin{equation}\label{statedependentapprox}
    \tau_{ij}(t) = \ep \big| q_i(t) - q_j(t) \big| + O(\ep^2).
  \end{equation}

 Keeping only the first order approximation in
 \eqref{statedependentapprox} makes \eqref{electrodynamics} an SDDE,
 but with \eqref{delaydefined}, the delay depends implicitly on the
 trajectory.
 
Note that it is not true that $\tau_{ij} = \tau_{ji} $ even if this
symmetry is true in the first order approximation
\eqref{statedependentapprox}.
\item
  In the case that $\tau_{ij} = 0$ and that the external forces are
  autonomous, the energy is conserved.  This has two consequences:
  \begin{itemize}
    \item
    In the autonomous case, the periodic orbits do not satisfy the
    hypothesis \eqref{H1}. Hence, we will only make precise statements
    in the case of time periodic external fields. In this case (very
    well studied in accelerator physics, plasma, etc.), there are many
    examples of periodic orbits satisfying assumption \eqref{H1''}, so
    that the results presented here are not vacuous.
\item 
  If the external potential and external magnetic fields are bounded,
  the periodic orbits of finite energy and away from collisions
  satisfy $|q_i(t) - q_j(t)| \ge \delta $, $i \ne j$, and $| q'_i(t) |
  \le c -\delta$. We will assume these two properties.
  \end{itemize} 
\end{enumerate}

Denoting $y(t) \coloneq (q_1(t),\ldots, q_N(t), q'_1(t),\ldots q'_N(t)
)$, we can write the equation \eqref{electrodynamics} in the form of
\eqref{smalldelay} with the delays being implicitly defined. Note that
there are $N(N-1)$ delays in total.

\begin{remark} 
Even if we formulate the result for the retarded potentials, we point
out that the mathematical treatment of Maxwell equations admits also
advanced potentials.

It is customary to take only the retarded potentials because of
``physical reasons'' which are relegated to footnotes in most
classical electrodynamics books. More detailed discussions appear in
\cite{Rohrlich, Spohn04}.  Note, that selecting only retarded
potentials breaks, even at the classical level, the time-reversibility
present in Maxwell's and Newton's equations. Mathematically any
combination of advanced and retarded potentials would make sense from
Maxwell equations.  Indeed, \cite{WheelerF45} proposes a theory with
half advanced and half retarded potentials.

We do not want to enter now into the physical arguments, which should
be decided by experiment (we are not aware of explicit experimentation
of these points). We just point out that the mathematical theory here
and the asymptotic expansions \cite{CasalCL20} applies to retarded,
advanced, or combination of advanced and retarded potentials.
\end{remark}

\subsection{Mathematical results for electrodynamics}
\label{sec:mathematicalresult}

In this section, we will collect the ideas we have been establishing
and formulate our main result for the model \eqref{electrodynamics}.
Note that we formulate the result only for periodic external fields,
since when the external fields are time-independent, energy is
conserved which prevents periodic orbits from satisfying assumption
\eqref{H1}.

We will assume that there exists $0<\xi_1 <1$, and $\xi_2>0$, such
that for all $t$:
\begin{equation} 
\label{trajectory_condition}
\begin{split} 
& |q'_j(t)| \le \xi_1 c  \\ 
& |q_i(t) - q_j(t) | \ge \xi_2
\end{split} 
\end{equation}

Note that \eqref{trajectory_condition} implies that the internal
forces and the masses are analytic around the trajectory.  Therefore,
the regularity assumptions for the equation concern only the external
fields.

\begin{theorem}\label{thm:electrodynamics} 
Denote $\ep = 1/c$. Consider the model \eqref{electrodynamics} with
the delays defined in \eqref{delaydefined}.

Assume that for $\ep = 0$, the resulting time periodic ODE has
periodic solution satisfying hypothesis \eqref{H1''} as well as
\eqref{trajectory_condition}. Assume that the external fields
$A_{\text{ext}}$ are $C^{\ell + \Lip}$.

Then, for small enough $\ep$, we can find a $C^{\ell + \Lip} $
periodic solution of \eqref{electrodynamics}.

In case that the external fields are jointly $C^{\ell + \Lip}$ in
time, position, velocity, and in a parameter $\gamma$, the periodic
solutions are jointly $C^{\ell + \Lip}$ as functions of the variable
of the parameterization and the parameter $\gamma$.
\end{theorem} 

The proof follows the steps of Theorem~\ref{thm:smalldelay} once we
have the estimates on delays \eqref{rbounds} and \eqref{c0for_r},
which will be discussed in the next section.

\subsection{Some preliminary results on the regularity of the delay}

In this section, we study \eqref{delaydefined} as an equation for
$\tau_{ij}(t)$ when we prescribe the trajectories $q_i$ and $q_j$.
This makes precise the notion that the delay is a functional of the
whole trajectory.

In the following proposition, we collect the proofs of estimates that
establish \eqref{rbounds} and \eqref{c0for_r}.  Both follow rather
straightforwardly from considering \eqref{delaydefined} as a
contraction mapping.

\begin{proposition} \label{prop:tauestimates}
Let $q_i$ and $q_j$ be continuously differentiable trajectories that
satisfy \eqref{trajectory_condition}.

Then, for each $t\in \R$, we can find a unique $\tau_{ij}(t) > 0 $
solving \eqref{delaydefined}.

Moreover: 

If the trajectories $q_i$ and $q_j$ are $C^{\ell+\Lip}$, then the
$\tau_{ij}$ is $C^{\ell + \Lip}$.  There is an explicit expression
\begin{equation} 
\label{taurbounds} 
\| \tau_{ij} \|_{C^{\ell + \Lip}} 
\le \phi( \| q_i\|_{C^{\ell + \Lip}},  \| q_j\|_{C^{\ell + \Lip}}, \xi_1, \xi_2).
\end{equation}

Let $q_i, q_j, \bar q_i$, and $\bar q_j$ be trajectories satisfying
\eqref{trajectory_condition}.  Denote by $\tau_{ij}$ and $\bar
\tau_{ij}$ the solutions of \eqref{delaydefined} corresponding to
$q_i, q_j$ and to $\bar q_i, \bar q_j$, respectively.  Then we have:
\begin{equation} 
\label{tauc0} 
\| \tau_{ij} - \bar \tau_{ij} \|_{C^0} \le
C(\xi_1, \xi_2) \left( \| q_i - \bar q_i \|_{C^0} + \| q_j - \bar q_j\|_{C^0}\right).
\end{equation} 
\end{proposition} 

\begin{proof} 
Fix $t$ and, hence, $q_i(t) $. 
We treat \eqref{delaydefined} as a fixed point problem for the -- long
named -- unknown $\tau_{ij}(t) $ with the functions $q_i, q_j$ as well
as the number $t$ fixed.

The first part of the asumption \eqref{trajectory_condition} implies
that the RHS of \eqref{delaydefined}, as a function of $\tau_{ij}(t)$
has derivative with modulus bounded by $\xi_1< 1$. Hence, we can apply
the contraction mapping principle.  This establishes existence and
uniqueness.

Moreover, we can apply the implicit function theorem and obtain that
$\tau_{ij}(t)$ is as differentiable on $t$ as the RHS of
\eqref{delaydefined}. Furthermore, we can get expressions for
$\frac{d^k}{dt^k } \tau_{ij}(t) $ which are algebraic expressions
involving derivatives with respect to $t$ of $q_i(t),~q_j(t)$ up to
order $k$, and derivatives of $\tau_{ij}(t)$ up to order $k-1$. The
exact combinatorial formulas are very well known.  Using recurrence in
the order of derivatives, we obtain \eqref{taurbounds}.

To prove \eqref{tauc0}, we observe that since the contraction we used
before is uniform in $t$, we can consider the RHS of
\eqref{delaydefined} as a contraction in $C^0$.

We evaluate the RHS of \eqref{delaydefined} corresponding to $\bar
q_i$ and $\bar q_j$ on $\tau_{ij}$, note that
\begin{align*}
\bar q_i (\cdot) - \bar q_j( \cdot - \tau_{ij}(\cdot) ) = &
\left(\bar q_i(\cdot) - q_i(\cdot) \right) 
+ \left( q_j(\cdot -\tau_{ij}(\cdot) ) - \bar q_j( \cdot - \tau_{ij}(\cdot)\right) \\
&+\left( q_i (\cdot) -  q_j( \cdot - \tau_{ij}(\cdot))  \right).
\end{align*}

Hence, 
\[
\left\| \frac{1}{c}  \big|\bar q_i (\cdot) - \bar q_j( \cdot - \tau_{ij}(\cdot) )\big|-
\tau_{ij}(\cdot) \right\|_{C^0} 
\le
\frac{1}{c}\| q_i - \bar q_i \|_{C^0} + 
\frac{1}{c}\| q_j - \bar q_j \|_{C^0}
\]
From this, \eqref{tauc0} follows from the Banach contraction mapping. 
\end{proof} 

\begin{remark}
Notice that the delays $\tau_{ij}$'s contain small factor
$\frac{1}{c}$, so are the right hands of the inequalities
\eqref{taurbounds} and \eqref{tauc0}, as we can see in the proof
above. We can view $\tau_{ij}\coloneq\frac{1}{c} r_{ij}$ to fit in the
case of small delays.
\end{remark}

\section{The case of Hyperbolic periodic orbits}
\label{sec:hyperbolic} 

Our main result Theorems \ref{thm:per} and \ref{thm:smooth} are based
on the assumption \eqref{H1}, which is automatically satisfied when
the periodic orbit of the unperturbed equation is hyperbolic. Hence,
the main results of this section can be viewed as corollaries of
Theorems \ref{thm:per} and \ref{thm:smooth}. In fact, we need slightly
stronger assumptions in the regularity in this section.

In this section, we will introduce an operator, see
\eqref{inveqseparated}, which is slightly different from the one
introduced in Section \ref{sec:operator}.

Even if the operator considered in this section requires more
regularity in the finite dimensional case, it generalizes our results
to perturbations of partial differential equations, see Section
\ref{sec:pde}, to perturbations of Delay Differential Equations, and
to other solutions that we will not discuss here (quasi-periodic,
normally hyperbolic manifolds). We also note that the corrections
needed in this section can be independent of the period. This makes it
possible to develop a theory of aperiodic hyperbolic sets. We hope to come back to this problem.

\subsection{Dynamical definition of hyperbolic periodic orbits} 
It is a standard notion that a periodic orbit of the ODE $\dot x =
f(x)$ is hyperbolic when the following strengthening of \eqref{H1}
holds.

With the same notation as in Section~\ref{sec:assump}, we say that a
periodic orbit is hyperbolic if:

\begin{enumerate}
\renewcommand{\theenumi}{H1.1}
\renewcommand{\labelenumi}{(\theenumi)}
\item \label{H1.1} $\Phi(\theta_0+1;\theta_0)$ has a simple eigenvalue
  1 whose eigenspace is generated by $DK_0(\th_0)$.  Moreover, all the
  other eigenvalues of $\Phi(\theta_0+1;\theta_0)$ have modulus
  different from $1$.
\end{enumerate}

The assumption \eqref{H1.1} is equivalent to the following
evolutionary formulation \eqref{H1.1'} in terms of invariant
decompositions.  In the finite dimensional case, this formulation is
easily obtained by taking the stable and unstable spaces of the
monodromy matrix and propagating them by the variational equations.
In the infinite dimensional cases, similar formulations are obtained
using semi-group theory under appropriate spectral assumptions.

\begin{enumerate}
\renewcommand{\theenumi}{H1.1'}
\renewcommand{\labelenumi}{(\theenumi)}
 \item \label{H1.1'} For every $\theta \in \T$ there is a
   decomposition
\begin{equation} 
\label{decomposition} 
\R^n = E^s_\theta \oplus E^u_\theta \oplus E^c_\th, \quad E^c_\th=\Span \{ D K_0(\theta)  \},
\end{equation} 
depending continuously on $\theta$ such that $E^s _\theta$ is forward
invariant, $E^u _\theta$ is backward invariant under the variational
equation. Moreover, the forward semiflow (resp. backward semiflow) of
the variational equation is contractive on $E^s _\theta$ (resp. $E^u
_\theta$).
\end{enumerate}
\medskip

More explicitly, we can find families of linear operators 
\begin{equation*}
\begin{aligned} 
&\{ U^s_{\theta}(t)  \}_{\theta \in \T,t \in \R_+ }, &
U^s_{\theta}(t)&\colon E^s_\theta \rightarrow E^s_{\th+\omega_0 t}\quad t \in \R_+, \\ 
&\{ U^u_{\theta}(t)   \}_{\theta \in \T,t \in \R_- },& U^u_{\theta}(t)&\colon E^u_\theta \rightarrow E^u_{\th+\omega_0 t} 
\quad t \in \R_-,
\end{aligned} 
\end{equation*}
satisfying for all $\th\in\T$
\begin{equation} 
\begin{split}  
\partial_t  U^\sigma_{\theta}(t) &= 
Df(\k0(\om_0 t + \theta))   U^\sigma_{\theta}(t)  \quad \sigma \in \{s, u\}  \\
U^\sigma_{\theta}(0) &= \text{Id}\big|_{E^\sigma_\theta},
\end{split} 
\end{equation} 
and
\begin{equation}  \label{cocyle} 
U^\sigma_{\theta}(t + \tau)  =   U^\sigma_{\omega_0 t + \theta} ( \tau).
\end{equation} 
Moreover, there exist $C > 0$, $\mu_s >  0$, $\mu_u >  0$ such that
\begin{equation} \label{contraction} 
\begin{split} 
\| U^s_{\theta}(t) \| &\le C e^{-\mu_s t}  \phantom{,,}\quad t \ge 0, \\
\| U^u_{\theta}(t)\| &\le C e^{-\mu_u |t| } \phantom{,} \quad t \le  0. \\
\end{split} 
\end{equation} 

We can also define an evolution operator $U^c_{\theta}(t)$ in the
$E^c$ direction. Note that $U^c_{\theta}(1) =
\text{Id}|_{E^c_{\theta}}$.

\subsection{Main Result in Hyperbolic Case}\label{sec:mainhyper}

The first result in this case is that Theorem \ref{thm:per} is true if
assumption \eqref{H1} is changed to assumption \eqref{H1.1}(or
\eqref{H1.1'}), and assumption \eqref{H2.1} is strengthened to
\eqref{H2.1.1} as follows:
\begin{enumerate}
\renewcommand{\theenumi}{H2.1.1}
\renewcommand{\labelenumi}{(\theenumi)}
\item \label{H2.1.1} If $K\in U_\rho$ and $\om\in B_\delta$, then
  $\P(K, \om, \ga, \cdot)\colon \T\to \R^n$ is $C^{\ell+\Lip}$, with
  $\|\P(K, \om, \ga, \cdot)\|_{C^{\ell+\Lip}}\leq \phi_{\rho,\delta}$,
  where $\phi_{\rho,\delta}$ is a positive constant.
\end{enumerate}
Recall that $U_\rho$ is the ball of radius $\rho$ in the space
$C^{\ell+\Lip}(\T, \R^n)$ centered at $\k0$, and $B_\delta$ is the
interval in $\R$ with radius $\delta$ centered at $\om_0$.

The second result is that the results in Theorem \ref{thm:smooth} is
true if assumption \eqref{H1} is substituted by assumption
\eqref{H1.1}(or \eqref{H1.1'}), and assumption \eqref{H2.2} is
strengthened to (H2.2.1) as follows:
\begin{enumerate}
\renewcommand{\theenumi}{H2.2.1}
\renewcommand{\labelenumi}{(\theenumi)}
\item \label{H2.2.1} If $K\in \pazocal{U}_\rho$ and $\om\in
  \pazocal{B}_\delta$, then $\P(K, \om, \cdot, \cdot)\colon \T\times
  O\to \R^n$ is $C^{\ell+\Lip}$, with $\|\P(K, \om, \cdot,
  \cdot)\|_{C^{\ell+\Lip}}\leq \phi_{\rho,\delta}$, where
  $\phi_{\rho,\delta}$ is a positive constant.
\end{enumerate}
Recall that $\pazocal{U}_\rho$ is the ball of radius $\rho$ in the
space $C^{\ell+\Lip}(\T\times O, \R^n)$ centered at $\k0$, and
$\pazocal{B}_\delta$ is the ball in $C^{\ell+\Lip}(O, \R)$ with radius
$\delta$ centered at constant function $\om_0$.
\begin{remark}
We emphasize that the results in this section are weaker than Theorem
\ref{thm:per} and Theorem \ref{thm:smooth}, however, we want to
introduce a different operator in the proof which has applications in
ill-posed PDEs, see Section \ref{sec:pde}. Modification of the
operator will be useful in the study of other dynamical objects.
\end{remark}

\subsection{Proof}\label{sec:hyperp}
We proceed as in Section~\ref{sec:operator} and manipulate
\eqref{inveq} as a fixed point problem taking advantage of the
geometric structures assumed in $\eqref{H1.1'}$.

Given the decomposition as \eqref{decomposition},we define projections
$\Pi^s_\theta, \Pi^u_\theta, \Pi^c_\theta$ over the spaces
$E^s_\theta, E^u_\theta, E^c_\theta$.  We also use the notation
\[
 \kh^\sigma(\theta) \coloneq \Pi^\sigma_\theta \kh(\theta), \qquad
 \sigma \in \{s, u\}.
\]
Taking projections along the spaces of the decomposition, using the
variation of parameters formula, and taking the initial conditions to
infinity (this procedure is standard since \cite{Perron29}), we see
that \eqref{inveq} implies
\begin{equation}\label{inveqseparated}
\begin{split}
\omh&=\om_0\frac{\langle\Pi^c_{\th_0}\int^{\frac{1}{\om_0}}_{0} U^c_{\th_0+\om_0 t}(\frac{1}{\om_0}-t)B^\ep(\kh,\omh,\ga,\th_0+\om_0 t)dt, D\k0(\th_0)\rangle}{\left| D\k0(\th_0)\right|^2},\\
\kh^s(\theta) &= \int_{-\infty}^0   U^s_{\theta +  \omega_0 t} (-t)
\Pi^s_{\theta + \omega_0 t} B^\ep(\kh,\omh,\ga,\th_0+\om_0 t)  \, dt,  \\
\kh^u(\theta) &= -\int_{0}^{\infty}  U^u_{\theta +  \omega_0 t}( -t)
\Pi^u_{\theta + \omega_0 t}B^\ep(\kh,\omh,\ga,\th_0+\om_0 t)  \, dt.  \\
\end{split}
\end{equation}

Define the right hand side of \eqref{inveqseparated} as an operator of
$(\omh,\kh^s,\kh^u)$, one can get lemmas which are similar to Lemmas
\ref{lem:prop} and \ref{lem:contract}. Hence we can get a fixed point
of the operator in this case.

When the solutions of \eqref{inveqseparated} are smooth enough and
decay fast enough that we can take derivatives inside of the integral
sign (which will be the case of the fixed points that we produce), it
is possible to show, taking derivatives of both sides of
\eqref{inveqseparated} and reversing the algebra that the well behaved
fixed points of \eqref{inveqseparated} indeed are solutions of
\eqref{inveq}.

The remarkable aspect of \eqref{inveqseparated} is that we only need
$U^s$ for positive times, and $U^u$ for negative times. Hence, the
assumed bounds \eqref{contraction} imply that the indefinite integrals
in \eqref{inveqseparated} converge uniformly in the $C^{\ell + \Lip}$
sense. At the same time, we pay the price of requiring one more
derivative of $\P$ while using this operator.

Another important feature of the operator \eqref{inveqseparated} is
that it does not require many assumptions on the long term evolution
of the solutions (in Section~\ref{sec:operator} we use heavily that
the solutions we seek are periodic).  This makes it possible to use
analogues of \eqref{inveqseparated} in several other problems. We hope to come back to these questions in the near future.

\section{ The case of evolutionary equations with delays }
\label{sec:pde} 
\def \F{\mathcal{F}}
\def\integer{\mathbb{Z}}

In this section we extend the results on ODEs in the previous sections
to PDEs and other evolutionary equations (e.g. equations involving
fractional operators or integral operators).

The key observation is that, the previous treatments of periodic
solutions do not use much that the functions we are seeking take
values in a finite dimensional space. For example, the
Lemma~\ref{closure} is valid for functions taking values in Banach
spaces.  Hence, we will show that the methods developed in the
previous sections can be applied without much change to a wide class
of PDEs.

Indeed, since one of the points of the previous theory was to avoid
the discussion of the evolutions, the theory applies easily to PDEs
using only very simple results on the evolution of the PDE.

\begin{remark}\label{periodicPDE} 
In this paper, we will not discuss the existence of periodic solutions
of evolutionary equations before adding the delays.  There is already
a large literature in this area.

We point, however that in studying the periodic solutions of a PDE
(which lie in an infinite dimensional space), it is natural to
consider the periodic solutions of a finite dimensional truncation
(e.g. a Galerkin approximation).  The problem of going from the
periodic solutions of a finite dimensional problem to the periodic
solutions in an infinite dimensional space, has some similarity with
the problems dealt with in the first parts of this paper.

A framework that systematizes the passing from periodic solutions of
the Galerkin approximations to periodic solutions of the PDEs is in
\cite{FiguerasGLL}. The methods of \cite{FiguerasGLL} have some points
in common with the methods used in this paper. It bypasses the study
of evolutionary equations and just studies the functional equations
satisfied by a parametrization of a periodic orbit.  The methods in
\cite{FiguerasGLL} lead to computer-assisted proofs that have been
implemented in \cite{GameiroL, FiguerasL}.  Since the methods of
\cite{FiguerasGLL} and this paper have points in common, one can hope
to combine them and go from a periodic solution of Galerkin truncation
of the PDE to a periodic solution of the delay perturbation of the
PDE.
\end{remark}

\subsection{Formulation of the problem and preliminary results} 
We use the standard set up of evolutionary equations (See
\cite{diplodocus, Showalter97} ).

Consider problem of the form
\begin{equation} 
\label{PDEabstract}
\partial_t  u(t)  =  \F(u(t))  + \ep P(u(t),  u_t;\gamma),
\end{equation}
where $u(t)$, is the unknown and lies in a space $X$ consisting of
functions on a domain $\Omega$. The points in $\Omega$ will be given
the coordinate $x$, so that we can also consider $u(t,x)$ as a
function on $\R \times \Omega$.

The function space $X$ encodes regularity properties of the functions
as well as boundary conditions. In particular, changing the boundary
conditions, changes the space $X$ and therefore, the functional
analysis properties (e.g. spectra) of the operators acting on it.

The operator $\F$ is a (possibly nonlinear) differential (or
fractional differential etc.)  operator.

As before (and contrary to the standard use in PDEs where $u_t$
denotes partial derivative), we use $u_t$ to denote a segment of the
solution, which can be related with history or future.  For $s \in
[-h,h]$, $u_t(s) = u(t + s) $, so that $u_t \in \mathcal{R}([-h,h],
X)$, a space of regular functions on $[-h,h]$ with values in $X$.  To
denote derivatives with respect to time we will always use $\partial_t
u$.

We consider $P\colon X\times  \mathcal{R}([-h,h], X)\times \R\to X$.

It  is  useful to  think heuristically of 
\begin{equation}\label{PDEevolution} 
\partial_t u = \F(u)
\end{equation}
as a differential equation in $X$ and indeed, our results will be
based on this heuristic principle.  To make sense of this heuristic
principle we have to overcome the problem that in the interesting
applications (See e.g. Section~\ref{sec:examples}), $\F$ is highly
discontinuous (involving derivatives) and not defined everywhere so
that the standard tools for smooth ODEs do not apply, but this is a
well studied problem.

A research program which became specially prominent in the 60's shows
that one can recover many of the results (existence, dependence on
initial conditions, etc.)  for the equation \eqref{PDEevolution} by
assuming functional analysis properties of the operator $\F$, see
\cite{diplodocus, BersJS64,Showalter94, Showalter97, Henry81, SellY02,
  Chueshov02}.  Of course, the verification of the functional analysis
assumptions in concrete examples, requires some hard analysis. One of
the subtle points of this program is that the notion of solutions may
be redefined to be weak or mild solutions.

Even if we will use the language and some material from the above
program, we will take a different point of view.
\begin{itemize} 
\item  
In this paper, we will not be interested in the theory of existence
and well-posedness for {\bf ALL} the possible initial conditions.
\item 
Indeed, because we are not going to discuss the initial value
problems, we can consider situations where the set of initial
conditions for the delay problems are not clear. Nevertheless, we can
get existence of smooth solutions.
\item
Since we are only aiming to produce some particular solutions, one
gets stronger results by taking more reduced spaces so that the
solutions are more regular and can be understood in the classical
sense.  In particular, in all the cases we will consider, the
functions and their derivatives will be bounded. (This happens,
e.g. if $X$ is a Sobolev space of high enough order.)

This is in contrast with the general theory of existence and
uniqueness, where the figure of merit is considering a more general
space of initial conditions.
\item
A more elaborate set-up for existence of evolutions that includes also
functional differential equations is \cite{Wu96}. In this paper,
however, we will avoid discussing the evolution of the Functional
differential equations and need only some results on the evolution of
the PDE.
\end{itemize}

\subsection{Overview of the method} 
Roughly, we will formulate analogues of the operator $\gaep$ in
\eqref{op1} and \eqref{op2} as well as the operator in
\eqref{inveqseparated} and verify that similar contraction argument
can be carried out.

The requirements of the above program on the theory of existence are
very mild.  The operator $\gaep$ only requires the existence of
solutions of the variational equation for finite time.  The operators
formulated in \eqref{inveqseparated} only require the existence of
partial evolutions (forward and backward evolutions in complementary
spaces), which allows to consider ill-posed equations, see
Section~\ref{sec:illposed}. Moreover, the smoothness requirements on
the delay terms are very mild.

\subsection{Examples}
\label{sec:examples}
In this section, we will present some examples which are
representative of the results we establish and which have appeared in
applications.

Even if we hope that this section can serve as motivation, from the
purely logical point of view, it can be skipped.  Of course, our
results apply to many more models and this section is not meant to be
an exhaustive list but to provide some intuition.

\subsubsection{Delay Perturbations}
One example of delay perturbation which considers long range
interaction is
\begin{equation}\label{longrange}
P(u(t), u_t;\gamma)=\int_{\R^d} K(x,y)\cdot u(t,x) \cdot u(t
-\frac{1}{c} |x -y|, y) dy.
\end{equation}
 
This models a situation in which the position $x$ interacts with
position $y$ with a strength $K(x,y)$, with the interaction taking
some time (proportional to the distance) to propagate. In
\eqref{longrange} we have denoted by $c$ the speed of propagation of
the signal, which is assumed to be constant.

Note that the interaction term could be more general than quadratic,
and may involve higher spacial derivatives thanks to the smoothing
property of solutions. Meanwhile, the speed of propagation of the
signal may not be constant (the propagation of signals may depend on
their strength).

Another example
\begin{equation}
\label{nonlocal}
P(u(t),  u_t;\gamma)= \int_0^\infty  G(s, u(t-s,x))\, ds,
\end{equation}
treating non-local interaction, is very typical in the modeling of
materials with memory effects (for example \emph{thixotropic}
materials) where the properties of the materials depend on the
history.  The effect of the previous state at present time often
decrease when the time delay grows. This is reflected on the function
$G(s, u)$ decreasing when $s$ (the delay in the effect) increases.

Of course, the mathematical theory that will be developed accommodates
more complicated effects such as $G$ depending on spatial derivatives
of $u$.

There are many other $P(u(t), u_t;\gamma)$ that we can consider. We only
need $P$ to satisfy some assumptions on regularity and Lipschitz
property, see (\ref{H2.1*}), (\ref{H3.1*}), and (\ref{H2.1.1*}), where we actually allow loss of regularity in the space variable.

In the coming sections, we see examples of unperturbed equations
\eqref{PDEevolution}.

\subsubsection{Parabolic equations} 
\label{sec:parabolic} 

Consider the equation for $u\colon \R\times \R^d\to \R$:
\begin{equation}\label{examplePDEparabolic}
\begin{split}  
&\partial_t u =  \Delta u +  N(x, u,\nabla u)\\
& u(t, x) = u(t, x +e)\quad \forall e \in \integer^d
\end{split}
\end{equation} 
with $N$ vanishing to quadratic order. For simplicity, we have imposed
periodic boundary conditions in space.

Notice that we have not imposed initial conditions at $t=0$ in example
\eqref{examplePDEparabolic}. Indeed, the initial conditions needed
require some thought.

As we will see, our treatment overcomes other possible complications
not mentioned explicitly so far.  We mention them because they are
natural in modeling and eliminating them from the literature may be
motivated by the need to have a more mathematically treatable problem.

Let us just mention briefly some small modifications.

\begin{itemize} 
\item The unknown $u$ could take values in $\R^d$. Note that
  considering systems rather than scalar equations makes a big
  difference in some PDE treatments (based on maximum principle), but
  it is not an issue in our case.
\item 
The papers \cite{KosovalicP19, KosovalicP19b} consider damped wave
equations perturbed by a delay.  From the functional analysis point of
view, the damped wave equations are similar to
\eqref{examplePDEparabolic}.
\end{itemize} 

\subsubsection{Kuramoto-Sivashinsky equations} 
\label{sec:KS} 
The model below is called the Kuramoto-Sivashinsky equation.
\begin{equation} \label{examplePDEKS}
\begin{split}  
&\partial_t u =  \Delta u +  \Delta^2 u + \mu \partial_x (u^2) \\
& u(t, x) = u(t, x +e)\quad \forall e \in \integer^d
\end{split} 
\end{equation} 

The Kuramoto-Sivashinsky equations appear as \emph{amplitude
  equations} for many problems arising in a variety of applications
(water waves, chemical reactions, interactive populations, etc.).

{From} the mathematical point of view, when $d =1$ (reduction of models with more
variables), the equation is
known to have an inertial manifold (all the solutions converge to a
finite dimensional manifold), which can be analyzed by finite
dimensional methods.  The equation \eqref{examplePDEKS} is known to
have many periodic solutions.  A very large number was identified by
non-rigorous, but reliable methods in \cite{Cvitanovic}.  Rigorous
periodic solutions have been established in many papers, including
bifurcations in \cite{ArioliK, Zgliczinski}. From the point of view of
this paper, it is interesting to note that \cite{FiguerasL, GameiroL}
use computer assisted proofs to establish the existence of periodic
orbits.

The equations discussed in the previous two sections are parabolic
PDEs so that indeed, the evolution is well defined and the solutions
gain smoothness. The linearized operator $\Phi$ that enters in
\eqref{op1} and \eqref{op2} is also smoothing.  Of course, for large
solutions, there could be finite time blow ups, but we are in the
regime of periodic solutions, which are well behaved.

\subsubsection{The Boussinesq equations in long wave approximation for water waves} 
In this section we present some physical equations that are ill-posed
in the sense that it is impossible to define an evolution for every
initial condition.  On the other hand, these equations may possess
many interesting and physically relevant solutions.

Since one of the main ideas of our treatment of FDEs is to bypass the
evolution, we obtain results on delay perturbations of ill-posed
equations.  This indeed highlights the difference of the present
method with the methods in evolution equations.

The material of this section is somewhat more sophisticated than the
rest of the paper and does not affect any of the other results.

Consider the equation for $u\colon \R\times \R\to \R$, derived in
\cite{Boussinesq72} as a long wave approximation for water waves.

\begin{equation} \label{Boussinesq} 
\partial_t^2 u  =  \mu \partial_x^4 u + \partial_x^2 u +  (u^2)_x    
\quad u(t,x +1) = u(t,x) 
\end{equation} 

This equation \eqref{Boussinesq} can be written as an evolution
equation of the form \eqref{PDEevolution} as follows:
\begin{equation}\label{Boussinesqevolution} 
\begin{split} 
\partial_t u  &=  v \\  
\partial_t v  &=  \mu \partial_x^4 u +\partial_x^2 u +  (u^2)_x  \\ 
& u(t,x+1) = u(t,x); \quad v(t,x+1) = v(t,x)  
\end{split} 
\end{equation} 
The linear part of the evolution is 
\begin{equation}
\label{Boussinesqevolutionlinear} 
\begin{split} 
\partial_t u  &=  v \\  
\partial_t v  &=  \mu \partial_x^4 u + \partial_x^2 u\end{split} 
\end{equation} 

Equations similar to \eqref{Boussinesq} have also appeared in other
contexts.  In water wave theory, $\mu > 0$, this leads to
\eqref{Boussinesq} being ill-posed. Indeed, consider the linear part
of the equation, the coefficient of the $k$-th Fourier mode $\hat u_k$
satisfies $\frac{d}{dt^2}\hat u_k = ( \mu k^4 - k^2)\hat u_k$, which
leads to exponentially growing solutions either in the future or in
the past.

Nevertheless, it is well known that the Boussinesq equation contains
many physically interesting solutions, including traveling waves and
other periodic and quasi-periodic solutions that are not traveling
waves. Notably, it contains a finite dimensional manifold (local
center manifold) which is locally invariant and on which solutions can
be defined till they leave the local center manifold \cite{Llave09,
  LlaveS19, ChengL20}. In particular, the periodic and quasi-periodic
solutions in the local center manifold are defined for all times.

For our purposes, the Boussinesq equation \eqref{Boussinesq} is
Hamiltonian, so that all the periodic solutions have monodromy with
eigenvalues $1$ -- corresponding to the conservation of the energy --
which make them unsuitable for the present version of our theory.
Hence, we will consider, for $u\colon \R\times \R^d\to \R$, mainly
time periodic perturbations of \eqref{Boussinesq}, which following the
notation in \cite{ChengL20}, we write as:
\begin{equation}\label{Boussinesq29}
\begin{aligned}
\left\{
\begin{array}{l}
\partial_t \theta= \omega \\
\partial_t^2 u =\mu \Delta^{2}u+\Delta u+N_{1} (\theta,x)+N_{2}(\theta,x)u+N_{3}(\theta,x,u,\nabla  u,\Delta u),\\
t\in\mathbb{R},\ \ \theta\in\mathbb{T},\quad x\in \mathbb{T}^d,
\end{array}
\right.
\end{aligned}
\end{equation}

The model \eqref{Boussinesq29} can be a long wave approximation of a
water wave model perturbed periodically.  These are physically
sensible long wave approximations of a water wave subject to periodic
forcing (e.g. waves in the ocean subject to tides or water waves in a
vibrating table -- Faraday experiment).

The result of \cite{ChengL20} implies, under very mild regularity
assumptions on $N_1, N_2, N_3$, that there is a finite dimensional
local center manifold of \eqref{Boussinesq29} which is locally
invariant.

This local center manifold is modeled on $\mathbb{T} \times
\mathbb{R}^n$.  The periodic solutions in the manifold are defined for
all time.  For specific forms of $N$, it is possible to prove the
existence of periodic orbits of \eqref{Boussinesq29}, which are
non-degenerate in the center manifold.

A natural space to consider \eqref{Boussinesqevolution} is $(u, v) \in
X\coloneq H^r \times H^{r -1} $ for sufficiently large $r$.  Even if
it is impossible to define an evolution of the linear part
\eqref{Boussinesqevolutionlinear} in the full space $X$, it is easy to
show using Fourier analysis that there are two complementary spaces in
which one can define the evolution forwards and backwards.  A
remarkable result in \cite{LlaveS19, ChengL20} is that this splitting
with partial evolution operators persists in the linearization near
periodic orbits, provided that they stay close to the origin.

\subsection{Result for well-posed PDE}
\label{sec:wellposed} 

The Theorem~\ref{thm:wellposed} will be our main result for well-posed
PDEs.  Essentially, the assumptions of the theorem are that we can
formulate the functional equation in \eqref{op1} and \eqref{op2} and
that the delay term prossesses enough regularity so that the argument
we used to prove Theorem~\ref{thm:per} goes through unchanged.

Therefore, the proof of Theorem~\ref{thm:wellposed} is a trivial
walk-through.  On the other hand, the fact that the assumptions are
satisfied in the cases \eqref{examplePDEparabolic},
\eqref{examplePDEKS} for some choices of spaces $X$ is not trivial and
will be discussed in Section~\ref{sec:verification}.  Of course,
similar verifications can be done in other models.

The only subtlety is that we will use the \emph{two spaces approach}
of \cite{Henry81}. (See also \cite{Taylor3,Chue} for a more
streamlined and refined version.)  This allows to consider
perturbations which are unbounded but of lower order than the
evolution operator.  For example in \eqref{examplePDEparabolic}, the
nonlinearity involves the first derivatives taking advantage of the
fact that the main evolution operator is of second order. In the case
of \eqref{examplePDEKS}, since the linear term is a fourth order
elliptic operator, the nonlinearity could involve terms of order up to
three.  As we will see, the two space approach also allows to lower
the regularity requirements of the delay term.  (See hypotheses in
Theorem~\ref{thm:wellposed}.)

\subsubsection{The two spaces approach} \label{sec:twospace}
The basic idea of the two spaces approach is that we study the
evolution equation using two spaces $X,Y$ consisting of functions with
different regularity.  In applications to PDE, often $X= H^{r +k}, Y =
H^{r}$ with $H^r$ the standard Sobolev spaces or the product of these
spaces.  In our case, we will take $r$ large enough so that the
solutions are classical, and the space $H^r$ enjoys properties that it
is a Banach algebra and the composition operator is smooth.

Differential operators, which are unbounded from a space to itself
become bounded from $X$ to $Y$.  Then, the main evolution operator,
smooths things out, such that it maps $Y$ to $X$ in a bounded way. Of
course, the bound of the evolution as an operator from the rough space
$Y$ to the smooth space $X$ depends on the time that the evolution has
been acting and becomes singular as the time goes to zero, but we
assume that there are bounds for the negative powers, which ensures
integrability.

\subsubsection{Setup of the result} 
\label{sec:setupPDEwellposed}

Consider the evolutionary PDE \eqref{PDEevolution}. Let $X,Y$ be
Banach spaces consisting of smooth enough functions satisfying the
boundary conditions imposed on \eqref{PDEevolution}. We will assume
that $Y$ consists of less smooth functions, such that $\F$ is a
differentiable map from space $X$ to space $Y$.  One consequence is
that $X$ has a compact embedding into $Y$.

Let $\k0\colon \T \to X$ be a parameterization of the periodic orbit
of \eqref{PDEevolution}. As in Section~\ref{sec:inveq}, we use the
notation $K(\theta)=\k0(\theta)+\kh(\theta)$ with $\kh\colon \T \to
X$, and we derive formally the equation \eqref{inveqPDE}.

\begin{equation}\label{inveqPDE}
\om_0 D\kh(\th)-D\F(\k0(\th))\kh(\th)=B^\ep(K,\om,\ga,\th)-\omh D\k0(\th),
\end{equation}
where
\begin{align}\label{Bdef2}
B^\ep(K,\om,\ga,\th)&\coloneq N(\th, \kh)+\ep\P(K,\om,\ga,\th)-\omh D\kh(\th),\\
N(\th, \kh)&\coloneq \F(\k0(\th)+\kh(\th))-\F(\k0(\th))-D\F(\k0(\th))\kh(\th).\nonumber
\end{align}

\subsubsection{Statement of the result}

We first formulate an abstract result, Theorem~\ref{thm:wellposed},
whose proof is almost identical to the proof of
Theorem~\ref{thm:per}. The deep result is to verify that the
hypotheses of Theorem~\ref{thm:wellposed} hold in examples of
interest.  In Section~\ref{sec:verification}, we show that the
examples in Section~\ref{sec:examples} verify the hypotheses. We leave
the verification in other models of interest to the readers.

\begin{theorem}\label{thm:wellposed} 
Assume that when $\ep=0$, the equation \eqref{PDEabstract} has a
periodic orbit which satisfies:

\begin{itemize} 
\item 
The linearized equation around the periodic orbit admits a
solution. That is, for any $\th_0\in\T$ and $\th_0<\th\in\T$, there is
an operator $\Phi(\theta; \theta_0) $ mapping from $Y$ to $X$ solving
\begin{equation} \label{semigroup}
\om_0\frac{d}{d\th} \Phi(\th; \theta_0) =   D \F( \k0(\th)) \Phi(\th; \theta_0);
\end{equation} 
\item
\begin{itemize} 
\item 
$ 1 \in Spec(\Phi(1;0), X) \text{ is a simple eigenvalue}.  $
\item

The spectral projection on $ Spec(\Phi(1;0), X) \setminus \{1\}$ in
$X$ is bounded.\end{itemize}

\item 
The family of operators $\Phi$ is smoothing in the sense that it satisfies

\begin{equation} \label{smoothingbounds} 
\| \Phi(t; \theta_0) \|_{Y,X}  \le C (t -\theta_0)^{-\alpha}  \quad  0<\alpha < 1,
\end{equation}
where $\|\cdot\|_{Y,X}$ is the norm of an operator mapping from $Y$ to
$X$, $C$ is a constant.
\end{itemize}

We also need the following two assumptions on the delay
perturbation. Let $\ell>0$ be an integer. Denote the ball of radius
$\rho$ in the space $C^{\ell+\Lip}(\T, X)$ centered at $\k0$ as
$\pazocal{U}_\rho$, and the interval in $\R$ 
centered at $\om_0$ with radius $\delta$ as $B_\delta$.

\begin{enumerate}
\setcounter{enumi}{1}
\renewcommand{\theenumi}{H\arabic{enumi}.1*}
\renewcommand{\labelenumi}{(\theenumi)}
\item\label{H2.1*} If $K\in \pazocal{U}_\rho$ and $\om\in B_\delta$,
  then $\P(K, \om, \ga, \cdot)\colon \T\to Y$ is $C^{\ell-1+\Lip}$,
  with $\|\P(K, \om, \ga, \cdot)\|_{C^{\ell-1+\Lip}(\T, Y)}\leq
  \phi_{\rho,\delta}$, where $\phi_{\rho,\delta}$ is a positive
  constant.

\item \label{H3.1*} For $K,~K'\in \pazocal{U}_\rho$, and $\om,~\om'
  \in B_\delta$, there exists constant $\alpha_{\rho,\delta}>0$, such
  that for all $\th\in\T$,
\begin{equation*}
\phantom{AAA}\|\P(K, \om, \ga, \th)-\P(K', \om', \ga, \th)\|_Y\leq \alpha_{\rho,\delta}\max\left\{|\om-\om'|,\|K-K'\|_{C^0(\T,X)}\right\}.
\end{equation*}
\end{enumerate}

Then, for small enough $\ep$, the equation \eqref{PDEevolution} has a
periodic orbit, which is parameterized by a $C^{\ell+\Lip}$ map
$K\colon \T\to X$. $K$ is close to $\k0$ in the sense of
$C^{\ell}(\T,X)$.
\end{theorem} 
The proof of Theorem~\ref{thm:wellposed} is very easy. It suffices to
observe that, thanks to the hypotheses of the theorem, the operator
$\gaep$, defined in the same way as before, sends a ball in the space
$\R\times C^{\ell + \Lip}(\T, X)$ to itself and that in this ball, $\gaep$ is a
contraction under the norm of $\R\times C^{0}(\T, X)$.  Then, we apply
Lemma~\ref{closure}.

Similar to before, one can get smooth dependence on parameters result.

\subsubsection{Some remarks}

\begin{remark}
The assumption that equation \eqref{semigroup} admits solutions with
the bounds in \eqref{smoothingbounds} is rather nontrivial and its
verification in concrete examples requires PDE techniques.
\end{remark}

\begin{remark}
Thanks to \eqref{smoothingbounds}, $\Phi( 1; 0)$ is bounded from $Y$
to $X$ and, hence compact from $Y$ to $Y$.  Therefore, the spectrum
away from zero is characterized by the existence of finite dimensional
eigenspaces.

However, for an operator $A$ acting on two spaces $X\subset Y$, there
is no relation of $Spec(A,X)$ and $Spec(A,Y)$ in general.
\end{remark}

\begin{remark}
In our case, for the operator $\Phi( 1; 0)$, its point spectrum in
space $X$ agrees with its point spectrum in space $Y$. This is not
hard to see from the eigenvector equation and the smoothing effect of
the operator $\Phi( 1; 0)$.
\end{remark}

\subsubsection{Verification of the assumptions of Theorem~\ref{thm:wellposed} in 
some examples}  \label{sec:verification}

For the parabolic equations \eqref{examplePDEparabolic} and
\eqref{examplePDEKS}, a very elegant formalism is developed in
\cite{Henry81}. The case \eqref{examplePDEKS} will be simpler than
\eqref{examplePDEparabolic} since the linearized operator being higher
order leads to stronger smoothing properties of the evolution.

The space $Y$ will be $H^r$, a Sobolev space of high enough order. We
emphasize once again that for our purposes, the results are stronger
if the space is more restrictive.

The semigroup theory tells us that we can solve the equation
\eqref{semigroup} and that the solution is smoothing in the sense that
\begin{equation} \label{smoothing} 
\|\Phi(\theta; \theta_0)\|_{H^r, ~H^{r+a} } \le C(a) |\theta - \theta_0|^{-a}.
\end{equation} 
 
\subsection{Result for ill-posed PDE}
\label{sec:illposed} 

In this section, we show how one can get existence of periodic
solutions for delay perturbations of ill-posed PDEs.

We just need to assume that the linearized equation admits partial
evolutions (one evolution forward in time and another one backward in
time) defined in complementary spaces. If these evolutions are
smoothing, the methods of Section~\ref{sec:hyperbolic} apply without
change.

Again the deeper part is to show that the concrete examples satisfy
the assumptions.  In the case of the periodically forced Boussinesq
equation \eqref{Boussinesq29} with a periodic solution which is
hyperbolic, we will show that the periodic solution persists under
delay perturbation.  The assumption that \eqref{Boussinesq29} has a
hyperbolic periodic orbit is a non-trivial -- but easily verifiable in
concrete models -- assumption.  We note that the time independent
Boussinesq equation \eqref{Boussinesq} does not have hyperbolic
periodic orbits due to energy conservation.  Our results require
delicate regularity properties of the periodic orbits, which are
verified for all the bounded small solutions in \cite{ChengL20}.

Since the partial evolutions involve smoothing properties, we still
use the two spaces approach summarized in Section~\ref{sec:twospace}.
We have used the same set up as \cite{ChengL20} to help the reader
check for the applications.

\begin{remark} 
When the non-linear terms $N$ in \eqref{Boussinesq29} are analytic,
the periodic orbits are analytic. As mentioned in Remark
\ref{rmk:bootana}, we do not expect that the periodic orbits of the
perturbed equations are analytic.  So, we follow \cite{ChengL20} and
deduce the regularity of the periodic orbits from the $C^r$ regularity
of the center manifold.
\end{remark}

\subsubsection{Abstract setup for the study of ill-posed equations} 
We will assume that there is a periodic solution of the evolution
equation \eqref{Boussinesq29}, which satisfies the following
Definition~\ref{spectralnondegeneracy}. Definition~\ref{spectralnondegeneracy}
can be verified for the linear part of \eqref{Boussinesq29}, and is
shown to be stable under perturbations (which can be unbounded) in
\cite{LlaveS19, ChengL20}.  (Related notions of splittings and their
stability using a different functional analysis set up appear also in
\cite{ChowL96,LatushkinMR96}. We have found that the two spaces
approach is more concrete and easier to adapt to the delay case.)

Definition~\ref{spectralnondegeneracy} is motivated by an analogue of
hyperbolicity for ill-posed equations.  We do not assume that the
linearized equations define an evolution such as $\Phi$, but we assume
that there are two evolutions (one in the future and one in the past)
defined in complementary spaces. This is enough to follow the set up
introduced in Section~\ref{sec:hyperbolic} and formulate a fixed point
equation for the periodic orbit of the perturbed equation.

Let us make some remarks about some subtle technical points. 

$\bullet$ We assume that when these evolutions are defined, they are
\emph{smoothing}. That is, they take functions of a certain degree of
differentiability (in $x$) and map them into functions with more
derivatives.  As shown in \cite{LlaveS19, ChengL20}, this allows to
show that these structures are stable under perturbations, which can
be unbounded but are of lower order.  This generality is important in
the treatment of examples such as \eqref{Boussinesq} since it allows
to show that the periodic solutions constructed in the above papers
satisfy Definition~\ref{spectralnondegeneracy}.

$\bullet$ It is important to note that
Definition~\ref{spectralnondegeneracy} only needs to be applied to the
periodic orbits of the problem without the delay.  In this section the
unperturbed problem will be a PDE, which is exactly the case discussed
in \cite{LlaveS19, ChengL20}.  As in Section~\ref{sec:hyperbolic}, the
invariant splitting will be used to set up a functional equation and
it will remain fixed, so that once we verify the existence in the
unperturbed case, it does not get updated.

$\bullet$ Both \cite{LlaveS19, ChengL20} consider situations more
general than periodic orbits.  The paper \cite{LlaveS19} considers
quasi-periodic orbits and \cite{ChengL20} considers bounded orbits. In
the case of quasi-periodic (in particular periodic) orbits, it is
natural in the examples considered to assume that the bundles are
analytic.  For orbits with a time-dependence more complicated than
periodic, it is natural to assume only finite regularity.  In this
paper we have adopted the definition in \cite{LlaveS19}, which
includes analyticity, since it applies to the examples we have in
mind. Notice, however that the solutions of the delay equation will
only be shown to be finitely differentiable and depend regularly on
parameters in finite differentiable topologies. Indeed, we do not
expect that the solutions of the delay problem will be analytic.  See
Remark~\ref{rmk:bootana}.

\begin{definition}\label{spectralnondegeneracy}
Let $X\subset Y$ be two Banach spaces. We say that an embedding
$K_0\colon \mathbb{T}_{\rho}\rightarrow X$ is spectrally nondegenerate
if for every $\theta$ in $\mathbb{T}$, we can find splittings:
\begin{equation}\label{splitting}
\begin{split}
X=X_{\theta}^{s}\oplus X_{\theta}^{c}\oplus X_{\theta}^{u}\\
Y=Y_{\theta}^{s}\oplus Y_{\theta}^{c}\oplus Y_{\theta}^{u}
\end{split}
\end{equation}
with associated bounded projections on $X$ and $Y$. (We will abuse the
notation and use $\Pi_{\theta}^{s,c,u}$ to denote the projections as
maps in ${L}(X,X)$ or in ${L}(Y,Y)$.)  The projections depend
analytically on $\theta \in \mathbb{T}_\rho$, and have continuous
extensions to the closure $\mathbb{T}_\rho$. Spaces
$X_{\theta}^{s,c,u}$ and $Y_{\theta}^{s,c,u}$ have the following
properties.

\begin{itemize}
\item We can find families of operators
\begin{equation*}
\begin{aligned}
 &U_{\theta}^{s}(t)\colon Y_{\theta}^{s}\rightarrow X_{\theta+\om_0 t}^{s},\quad t>0,\\
&U_{\theta}^{u}(t)\colon  Y_{\theta}^{u}\rightarrow X_{\theta+\om_0 t}^{u},\quad t<0,\\
&U_{\theta}^{c}(t)\colon  Y_{\theta}^{c}\rightarrow X_{\theta+\om_0 t}^{c},\quad t\in\mathbb{R}.
\end{aligned}
\end{equation*}
\item The operators $U_{\theta}^{s,c,u}(t)$
are cocycles over the rotation satisfying
\begin{equation}\label{cocycles}
\begin{aligned}
U_{\theta+\om_0 t}^{s,c,u}(\tau)U_{\theta}^{s,c,u}(t)=U_{\theta}^{s,c,u}(\tau+t).
\end{aligned}
\end{equation}

\item The operators $U_{\theta}^{s,c,u}(t)$ are smoothing in the time
  direction where they can be defined and they satisfy assumptions in
  the quantitative rates. There exist constants $\alpha_1,\alpha_2\in
  [0,1)$, $\beta_1, \beta_2, \beta_3^+, \beta_3^- > 0$ with $ \beta_1
    > \beta_3^-$, and $ \beta_2 > \beta_3^+$, and $C>1, $ independent
    of $\theta,$ such that the evolution operators satisfy the
    following rate conditions:
\begin{equation}\label{stablecocycle}
\begin{aligned}
\|U_{\theta}^{s}(t)\|_{\rho,Y,X}\leq C e^{-\beta_{1}t} t^{-\alpha_{1}},
\quad t>0,
\end{aligned}
\end{equation}
\begin{equation}\label{unstablecocycle}
\begin{aligned}
\|U_{\theta}^{u}(t)\|_{\rho,Y,X}\leq C e^{-\beta_{2}|t|} |t|^{-\alpha_{2}},
\quad t<0,
\end{aligned}
\end{equation}
and
\begin{equation}\label{centercocycle}
\begin{aligned}
&\|U_{\theta}^{c}(t)\|_{\rho,Y,X}\leq C e^{\beta_{3}^{+}t} ,\quad t\geq0,\\
&\|U_{\theta}^{c}(t)\|_{\rho,Y,X}\leq C e^{\beta_{3}^{-}|t|} ,\quad t\leq0.
\end{aligned}
\end{equation}

\item
The operators $U_{\theta}^{s,c,u}(t)$ are solutions of the variational
equations in the sense that
\begin{equation}\label{integralsolution}
\begin{aligned}
U_{\theta}^{s}(t)=Id+\int_{0}^{t}D\F^{s}(\k0(\theta+\om_0\tau))U_{\theta}^{s}(\tau)d\tau,\quad t>0,\\
U_{\theta}^{u}(t)=Id+\int_{0}^{t}D\F^{u}(\k0(\theta+\om_0\tau))U_{\theta}^{u}(\tau)d\tau, \quad t<0,\\
U_{\theta}^{c}(t)=Id+\int_{0}^{t}D\F^{c}(\k0(\theta+\om_0\tau))U_{\theta}^{c}(\tau)d\tau, \quad t\in \mathbb{R}.
\end{aligned}
\end{equation}
\end{itemize}

In this paper, we will also need:
\begin{itemize}
\item The space $X^c$ is unidimensional and it is spanned by the
  direction of the evolution along the periodic orbit.
\end{itemize}
\end{definition}

Recall that $\pazocal{U}_\rho\subset C^{\ell+\Lip}(\T, X)$ is the ball of radius
$\rho$ centered at $\k0$, and $B_\delta\subset\R$ is the interval centered at $\om_0$ with radius $\delta$. Compared with the hypothesis for well-posed equations in
(\ref{H2.1*}), we make similar but slightly stronger assumption on the
delay term:

\begin{enumerate}
\renewcommand{\theenumi}{H2.1.1*}
\renewcommand{\labelenumi}{(\theenumi)}
\item \label{H2.1.1*}If $K\in \pazocal{U}_\rho$
  and $\om\in B_\delta$, then
  $\P(K, \om, \ga, \cdot)\colon \T\to Y$ is $C^{\ell+\Lip}$, with
  $\|\P(K, \om, \ga, \cdot)\|_{C^{\ell+\Lip}(\T,Y)}\leq
  \phi_{\rho,\delta}$, where $\phi_{\rho,\delta}$ is a positive
  constant.
\end{enumerate}

\subsubsection{Statement of the result}

\begin{theorem} \label{thm:ill-posed} 
Assume that we have an evolution equation \eqref{PDEevolution} that
admits a periodic solution satisfying
Definition~\ref{spectralnondegeneracy}, and that we perturb by delay
terms satisfying assumptions (\ref{H2.1.1*}) and (\ref{H3.1*}).

Then, for sufficiently small $\ep$, the equation \eqref{PDEabstract}
has a periodic solution of frequency $\omega$, which is parameterized
by a $C^{\ell+\Lip}$ map $K\colon \T\to X$. $K$ is close to $\k0$ in
the sense of $C^{\ell}(\T,X)$.
\end{theorem}

The proof of Theorem~\ref{thm:ill-posed} follows the same line as in
Section \ref{sec:hyperp}. We work with the fixed point equation
\eqref{inveqseparated}. Using that we have evolution $U^c$ and partial
evolutions $U^s$ and $U^u$ for the linearized equations satisfying
Definition \ref{spectralnondegeneracy}, we can find solution to
equation \eqref{inveqseparated}, with $\omh\in \R$ and $\kh^s$,
$\kh^u\in C^{\ell+\Lip}(\T,X)$.

As we have discussed, the regularity properties are verified for
concrete examples in \cite{ChengL20} for the time-perturbed Boussinesq
equation~\eqref{Boussinesq29}.

\appendix
\section{Regularity properties} 
\label{app:regularity} 
One of the sources of complication in the study of delay equations --
especially state dependent delay equations -- is that the equations
involve compositions, which have many surprising properties.  In this
appendix we collect a few of them.  A systematic study of the
composition operator in $C^r$ spaces which are the most natural for
our problem is in \cite{LO99}.

\subsection{Function Spaces}
Let $\ell$ be a positive integer, let $X$ be a Banach space and
$U\subset X$ be a an open set.  For functions on $U$ taking values in
another Banach space $Y$, we can define derivatives \cite{Dieudonne,
  LoomisS} and Lipschitz and H\"older regularity of the derivatives.

We recall that the $j$ derivative is a $j$-multilinear function from
$X^{\otimes j}$ to $Y$ and that there is a natural norm for
multilinear functions (supremum of the norm of the values when the
arguments have norm $1$).

For function $F$ defined in  domain in a Banach 
space, we denote by 
\[
\Lip{F} =  \sup_{x,y \in U, x\ne y} \| F(x) - F(y) \|_Y/\| x - y \|_X 
\]

\begin{definition} \label{def:ClLip} 
We say that $K\colon U \rightarrow Y$ is in $C^{\ell+ \Lip}(U, Y)$
when $K$ has $\ell$ derivatives and the $\ell$ derivative is
Lipschitz.

We endow  $C^{\ell+ \Lip}(U, Y)$ with the norm:

\begin{equation}\label{normdefined} 
\|K\|_{C^{\ell+\Lip}}=\max\left\{\max_{k=0,1,\dotsc, \ell}\big\{\sup_x\|D^k K(x)\|\big\}, \Lip(D^\ell K)\right\}
\end{equation}
which makes $C^{\ell+\Lip}$ into a Banach space.
\end{definition}

A similar definition can be written when $U$ is a Riemannian
manifold. In this paper we will use the case that $U = \T$.

\begin{remark} 
We note that Definition~\ref{def:ClLip} assumes uniform bounds of the
derivatives in the whole domain. There are other very standard
definitions of differentiable sets that only assume continuity and
bounds in compact subsets of $U$. Even when $U = \R^n$ these
definitions (e.g. Whitney topology, very natural in differential
geometry) do not lead to $C^{\ell + \Lip}$ being a Banach space and we
will not use them.
\end{remark}

\subsection{ Simple estimates on Composition}
We will need the following property of the composition operator, one
can refer to \cite{LO99} for more details.
\begin{lemma}\label{lem:com}
Let $X,Y,Z$ be Banach spaces. 
Let $E \subset X$, $F \subset Y$ be open subsets.

Assume that: $g\in C^{\ell+\Lip}(E, Y)$, $f\in C^{\ell+ \Lip}(F,Z)$
and that $g(E) \subset F$ so that $f\circ g$ can be defined.  Then,
$f\circ g \in C^{\ell+ \Lip}(E,Z)$, and
\begin{equation} \label{eqn:comp} 
\|f\circ g\|_{C^{\ell+ \Lip}(E,Z)}\leq M_\ell \|f\|_{C^{\ell+ \Lip}(F,Z)}\big(1+\|g\|_{C^{\ell+ \Lip}(E,Y)}^{\ell+1}\big)
\end{equation}

\end{lemma}

The proof of Lemma~\ref{lem:com} just uses the Faa-Di-Bruno formula
for the derivatives of the composition. To control the Lipschitz
constant of the $\ell$ derivative, we use that the Lipschitz constant
of product and composition satisfy the same formulas as those of the
derivative with an inequality in place of equality.

In \eqref{eqn:comp} we can take any set $F$ that contains $g(E)$. The
results are sharper when we take $F$ as small as possible.

\subsection{The mean value theorem} 

\begin{definition} \label{def:compensated} 
We say that an open set $U \subset X$ is a compensated domain when it
is connected, and there is $C> 0 $ such that for any $x,y \in U$,
there is a $C^1$ path $\gamma \subset U$ such that
\[
\text{length}(\gamma)  \le C \| x - y\| 
\]
\end{definition} 
In particular, a convex domain is compensated with $C =1 $. 

We also recall the fundamental theorem of calculus. 
\begin{theorem}  \label{FTC} 
Assume that $U \subset X$ is open connected, $F\colon U \rightarrow X$
is a $C^1$ function, $x, y \in U$ and that $\gamma$ is a $C^1$ path
joining $x,y$.  Then
\[
F(x) - F(y) = \int_{0}^1 DF(\gamma(t))\, D\gamma(t)  \, dt 
\]
\end{theorem} 

As a corollary of Theorem~\ref{FTC} we have that
\[
\| F(x) - F(y) \| \le  \|DF\|_{C^0} \cdot\text{length}(\gamma)  \le  \|F\|_{C^1} \cdot \text{length}(\gamma)
\]
If the domain $U$ is compensated, we obtain that
\[
\| F(x) - F(y) \| \le  C  \|F\|_{C^1} \|x - y\|
\]
In particular, $C^1$ functions on compensated domains are Lipschitz.

The conclusion that $C^1$ implies Lipschitz, is not true if the domain
is not compensated. It is not difficult to obtain examples of domains
where $C^1$ functions are not continuous even when $X = \R^2$.

\begin{lemma}\label{compositionLipschitz}
Assume that for some $\ell \ge 1$, $ \| f \|_{C^{\ell+ \Lip}} \le A $,
$\|g_1\|_{C^{\ell -1 + \Lip}}$, $\|g_2\|_{C^{\ell - 1 + \Lip}} \leq
B$.  Then:
\begin{equation} \label{Lipschitbound} 
 \|f \circ  g_1 - f\circ g_2 \|_{C^{\ell -1 + \Lip}}\leq C(A,B) \| g_1 -  g_2 \|_{C^{\ell -1 + \Lip}}
\end{equation} 
\end{lemma} 

\begin{proof} 
By the fundamental theorem of calculus we have pointwise
\[
f\circ g_1  - f\circ g_2  = \int_0^1  Df(g_2 + t(g_1 -g_2) ) (g_1 - g_2) \, dt
\]
If we interpret the above as identity among functions we have 
\[
\| f\circ g_1  - f\circ g_2  \|_{C^{\ell -1 + \Lip}}\leq \int_0^1 C  \| Df(g_2 + t(g_1 -g_2) )\|_{C^{\ell -1 + \Lip}} \cdot \| (g_1 - g_2)\|_{C^{\ell -1+ \Lip}} \, dt
\]
Using Lemma~\ref{lem:com}, $ \|Df(g_2 + t(g_1 -g_2) )\|_ {C^{\ell -1 +
    \Lip}}$ is bounded by a function of $A$ and $B$, we are done.
\end{proof}

\subsection{Interpolation}
We quote the following result from \cite{Had98, Kol49}.  See
\cite{LO99} for a modern, very simple proof valid for functions on
compensated domains in Banach spaces.

\begin{lemma}\label{lem:inter}
Let $U$ be a convex and bounded open subset of a Banach space $E$, $F$
be a Banach space. Let $r$, $s$, $t$ be positive numbers, $0\leq
r<s<t$, and $\mu=\frac{t-s}{t-r}$. There is a constant $M_{r,t}$, such
that if $f\in C^t(U,F)$, then
\[
\|f\|_{C^s}\leq M_{r,t}\|f\|_{C^r}^{\mu}\|f\|_{C^t}^{1-\mu}.
\]
\end{lemma}

\subsection{Closure Properties of $C^{\ell + \Lip}$ ball}

We quote a very practical result which appears as Lemma 2.4 in
\cite{Lan}. (This paper is largely reproduced as a chapter in
\cite{MarsdenM}. See Lemma (2.5) on p. 39)

\begin{lemma}\label{closure} 
Let $U \subset X$ be a compensated domain. 

Denote by $\mathbf{B} $ a closed ball in $C^{\ell + \Lip}(U,Y)$.  Let
$\{ u_n \}_{n\in \mathbb{N}} \subset \mathbf{B}$ be such that $u_n$
converges pointwise weakly to $u$.  Then, $u \in \mathbf{B}$.

Furthermore, the derivatives of $u_n$ of order up to $\ell$ converge
weakly to the derivatives of $u$.
\end{lemma}

We note that the hypothesis of Lemma~\ref{closure} are easy to verify
in operators that involve composition.  The propagated bounds just
amount to proving that the size of derivatives of composition of two
functions can be estimated by the sizes of the derivatives.  of the
original functions. The contraction properties are done under the
assumption that the functions are smooth so that one can use the mean
value theorem.

A similar result to Lemma~\ref{closure} is the following, which
appears as Lemma 6.1.6 in \cite[p. 151]{Henry81}.
\begin{lemma}\label{lem:henry}  
Let $U \subset X$ be an open set. Denote by $\mathbf{B} $ a closed
ball in $C^{\ell + \Lip}(U,Y)$.  Let $\{ u_n \}_{n\in \mathbb{N}}
\subset \mathbf{B}$ be such that $u_n$ converges uniformly to $u$.
Then, $u \in \mathbf{B}$.

Furthermore the derivatives of $u_n$ of order up to $\ell$ converge
uniformly to the derivatives of $u$ away from the boundary of $U$.
\end{lemma} 

Both Lemma~\ref{closure} and Lemma~\ref{lem:henry} remain true when we
replace the spaces of $C^{\ell+ \Lip}$ functions by H\"older spaces.

\begin{remark} 
It is instructive to compare the proofs of Lemma~\ref{closure} and
Lemma~\ref{lem:henry} in their original references.

The proof of \cite{Lan} is based on considering restrictions to
lines. Then, one can apply Arzela-Ascoli theorem and extract
converging subsequences.  The assumption of a weak pointwise limit
ensures that the limit is unique.  The uniformity of the $C^{\ell +
  \Lip}$ norms of the functions ensures the existence of derivatives
and the convergence.

The proof of \cite{Henry81} goes along different lines. It shows that
there are bounds on the derivatives by the $C^0$ norms and the size of
the ball. An alternative argument is to use interpolation inequalities
in Lemma~\ref{lem:inter}, which provides uniform convergence of the
derivatives on $U$ (also near the boundary).
\end{remark} 

As a consequence of Lemma~\ref{closure}, we have the following version
of the contraction mapping.

\begin{lemma}\label{lem:contraction} 
With the same notation of Lemma~\ref{closure}.

Assume $\mathcal{T}\colon \mathbf{B} \rightarrow \mathbf{B}$ satisfies
that there exists $\kappa < 1$ such that
\[
\| \mathcal{T}(u) - \mathcal{T}(v) \|_{C^0} \le \kappa \|u - v
\|_{C^0} \quad \forall u,v \in \mathbf{B}
\]
Then, $\mathcal{T} $ has a unique fixed point $u^*$ in $\mathbf{B}$. 

For any $u \in \mathbf{B}$, and $0 \le j \le \ell$ 
\[
\| \mathcal{T}^n(u) - u^*\|_{C^{j+\Lip}} \le C \kappa^{n\frac{\ell
    -j}{\ell +1}} \|\mathcal{T}(u) - u \|_{C^0}^{\frac{\ell -j}{\ell
    +1}}
\]
where $C$ is a constant that depends on the radius of the ball
$\mathbf{B}$ and $j$.

Furthermore, 
\[
\| u - u^*\|_{C^{j+ \Lip}} \le C (1 - \kappa)^{-\frac{\ell -j}{\ell
    +1}} \|\mathcal{T}(u) - u \|_{C^0}^{\frac{\ell -j}{\ell +1}}
\]
\end{lemma}  

\begin{proof} 
When $X$ is finite dimensional (or just separable),
Lemma~\ref{closure} is a corollary of Ascoli-Arzela theorem.  For any
subsequence of $u_n$ we can extract a sub-subsequence that converges
in $C^\ell$ sense.  The limit of this sub-subsequence has to be
$u$. It follows that the $u_n$ converges to $u$ in $C^\ell$ sense. It
then follows that the $\ell$-derivative is Lipschitz.

If $X$ is infinite dimensional, one can repeat the above argument
restricting to lines.  The uniform regularity assumed on $u_n$
translates to uniform regularity of $u$ restricted to lines.

We refer to \cite{Lan}. Indeed \cite{Lan} only needs to assume that
the sequence converges weakly pointwise. The convergence properties
are only used to guarantee the uniqueness of the limit obtained
through compactness (The paper \cite{Lan} is written when the domain
$U$ is the whole space, but this is not used).

Once we have the closure property, the existence of the unique fixed
point is as in Banach contraction. We observe that for any $u \in
\mathbf{B}$,
\[
\| \mathcal{T}^{n+1} (u) - \mathcal{T}^n(u) \|_{C^0} \le \kappa^n \|
\mathcal{T}(u) - u \|_{C^0}
\]

Using the interpolation inequalities Lemma~\ref{lem:inter} and that
the $C^{\ell+ \Lip}$ norms of the iterates are bounded, we obtain
\begin{equation}\label{geometricinterpolated}
\| \mathcal{T}^{n+1} (u) - \mathcal{T}^n(u) \|_{C^{j + \Lip} } \le C
\kappa^ {n\frac{\ell -j}{\ell +1}} \|\mathcal{T}(u) - u
\|_{C^0}^{\frac{\ell -j}{\ell +1}}
\end{equation} 

From this one obtains that $\mathcal{T}^n(u) - u = \sum_{k=1}^n
(\mathcal{T}^k(u) - \mathcal{T}^{k-1}(u) )$ is an absolutely
convergent series in the $C^{j + \Lip}$ sense. Let $u^*$ be the fixed
point. Using \eqref{geometricinterpolated} to estimate the series, we
obtain:
\[
 \|u - u^*\|_{C^{j+\Lip}} \le C(1 - \kappa^{\frac{\ell -j}{\ell
     +1}})^{-1} \|\mathcal{T}(u) - u \|_{C^0}^{\frac{\ell -j}{\ell
     +1}}.
\]

On the other hand, from the standard Banach fixed point theory, we
obtain that $\|u - u^*\|_{C^0} \le (1 - \kappa)^{-1} \| \mathcal{T}(u)
- u\|_{C^0}$.  By Lemma~\ref{lem:inter} we obtain
\[
 \|u - u^*\|_{C^{j+\Lip}} \le C(1 - \kappa)^{-\frac{\ell -j}{\ell +1}}
 \|\mathcal{T}(u) - u \|_{C^0}^{\frac{\ell -j}{\ell +1}}.
\]
It is easy to see that this bound is better than the previously
obtained one summing the series.
\end{proof} 
\bibliographystyle{alpha}
\bibliography{ref3}
\end{document}